\documentclass[12pt,reqno]{amsart} %us uses letter paper

\usepackage{times}

\usepackage[
  margin=1in
]{geometry}

\usepackage[pdftex]{graphicx}

\usepackage{array,multirow,makecell}
\usepackage{multicol}
\usepackage{subcaption}

\usepackage{array,multirow,makecell}
\setcellgapes{1pt}
\makegapedcells
\newcolumntype{R}[1]{>{\raggedleft\arraybackslash }p{#1}}
\newcolumntype{L}[1]{>{\raggedright\arraybackslash }p{#1}}
\newcolumntype{C}[1]{>{\centering\arraybackslash }p{#1}}

\usepackage{tikz}
\usetikzlibrary{tikzmark}

\usepackage{amssymb,amsfonts,amsthm}
\usepackage{marginnote}
\usepackage{comment}
\usepackage{enumitem}
\usepackage{mathrsfs}
\usepackage{changepage}

\usepackage[colorlinks=true, pdfstartview=FitV, linkcolor=black, citecolor=blue, urlcolor=blue]{hyperref}

\def\Xint#1{\mathchoice
{\XXint\displaystyle\textstyle{#1}}%
{\XXint\textstyle\scriptstyle{#1}}%
{\XXint\scriptstyle\scriptscriptstyle{#1}}%
{\XXint\scriptscriptstyle%
\scriptscriptstyle{#1}}%
\!\int}
\def\XXint#1#2#3{{\setbox0=\hbox{$#1{#2#3}{%
\int}$ }
\vcenter{\hbox{$#2#3$ }}\kern-.6\wd0}}
\def\barint{\, \Xint -} % \, corrects the \! used in the definition
\def\bariint{\barint_{} \kern-.4em \barint}
\def\bariiint{\bariint_{} \kern-.4em \barint}
\renewcommand{\iint}{\int_{}\kern-.34em \int} %\ minor space between the integrals
\renewcommand{\iiint}{\iint_{}\kern-.34em \int} %\ minor space between the integrals

\DeclareMathAlphabet{\mathcal}{OMS}{cmsy}{m}{n}

\theoremstyle{plain}
\newtheorem{goala}{Objective}

\newtheorem{resulta}{Result}

\theoremstyle{definition}

\newcommand{\R}{\mathbb{R}}

\newcommand{\ep}{\varepsilon}

\definecolor{darkgreen}{rgb}{0,0.5,0}
\definecolor{darkblue}{rgb}{0,0,0.7}
\definecolor{darkred}{rgb}{0.9,0.1,0.1}
\definecolor{lightblue}{rgb}{0,0.51,1}

\numberwithin{equation}{section}
\setlist[enumerate]{leftmargin=*}

\title[From concentration to quantitative regularity]{From concentration to quantitative regularity: a short survey of recent developments for the 
Navier-Stokes equations}
\author[T. Barker]{Tobias Barker}
\address[T. Barker]{Department of Mathematical Sciences, University of Bath, Bath BA2 7AY, UK}
\email{tobiasbarker5@gmail.com}

\author[C. Prange]{Christophe Prange}
\address[C. Prange]{Cergy Paris Universit\'e, Laboratoire de Math\'ematiques AGM, UMR CNRS 8088, France}
\email{christophe.prange@cyu.fr}
\date{\today}

\begin{document}

\maketitle

\begin{center}
\emph{In honor of Carlos Kenig's 70th birthday}
\end{center}

\begin{abstract}
In this short survey paper, we focus on some new developments in the study of the re\-gularity or potential singularity formation for solutions of the 3D Navier-Stokes equations. Some of the motivating questions are: Are certain norms accumulating/concentrating on small scales near potential blow-up times? At what speed do certain scale-invariant norms blow-up? Can one prove explicit quantitative regularity estimates? Can one break the criticality barrier, even slightly? We emphasize that these questions are closely linked together. Many recent advances for the Navier-Stokes equations are directly inspired by results and methods from the field of nonlinear dispersive equations.
\end{abstract}


\vspace{0.3cm}


\tableofcontents

\section{Introduction}

This short survey paper is concerned with recent developments in the study of the regularity or potential singularity formation for solutions of the 3D Navier-Stokes equations
\begin{equation*}
\partial_tU-\Delta U+U\cdot\nabla U+\nabla P=0,\qquad \nabla\cdot U=0.
\end{equation*}
We mostly focus on the whole-space case $\R^3$ or localize away from physical boundaries. We will mainly concentrate on two topics, concentration of solutions on small scales on the one hand and quantitative regularity on the other hand, and show how these subjects are related. The study of these questions is recent for the Navier-Stokes equations, all the main results in this paper were published in the last five years. 

In 2003 Escauriaza, Seregin and \v{S}ver\'{a}k \cite{ESS2003} were able to the blow-up of the critical borderline norm $L^3$ for potentially singular solutions of 3D Navier-Stokes. This paper was followed by a tremendous amount of works proving analogous results for many kind of evolution equations with a scaling symmetry\footnote{Divergence of critical norms near maximal time of existence for PDEs with a scaling symmetry does not follow from local well-posedness theory, see \cite{MR06}.} \cite{MR06,KM10,KM11,KV11,DY18,MizoguchiSouplet19,MT22,camliyurt2022scattering}\ldots\  These results are all qualitative, with the exception of the work by Merle and Rapha\"el \cite{MR06} for nonlinear Schr\"odinger, which gives a quantitative blow-up rate for a critical norm.

A breakthrough for the Navier-Stokes equations was achieved by Tao in 2019 \cite{Tao19}. He was able to explicitly quantify the rate of blow-up of the critical $L^3$ norm near a potential first-time singularity
\begin{equation}\label{e.quantbluprateTao}
\limsup_{t\rightarrow T^*}\frac{\|U(\cdot,t)\|_{L^3(\R^3)}}{\big(\log\log\log(T^*-t)\big)^{c}}=\infty,
\end{equation}
for a universal constant $c\in(0,\infty)$. Previously, only abstract quantitative results were known, which were based on abstract quantification of the seminal (qualitative) result of Escauriaza, Seregin and \v{S}ver\'{a}k \cite{ESS2003} and the use of persistence of singularities.

Some techniques described in this paper are directly inspired by methods introduced for nonlinear dispersive equations. In this vein, let us mention the `stacking of scales scheme' used to prove quantitative regularity estimates in Section \ref{sec.genstrat} and Section \ref{sec.quanttypeI} inspired by \cite{MR06}, and the `mild criticality breaking' result in Section \ref{sec.mildbreak} inspired by \cite{Bulut20}.

\subsection*{Explicit quantitative estimates: for what purpose?}

Blow-up rates and quantitative regularity estimates are two sides of the same coin. Let us outline four motivations for the study of quantitative regularity and blow-up rates:
\begin{enumerate}[label=(\arabic*)]
\item In the field of PDEs, there seems to be very few quantitative rates for critical norms near the maximal time of existence.
\item Quantitative regularity estimates with explicit bounds under a priori boundedness of critical norms enable to break (to some limited extent) the criticality barrier. For instance in Section \ref{sec.mildbreak}, we use the result of Tao \cite{Tao19} to derive a new regularity criteria in terms of a slightly supercritical Orlicz norm.
\item Blow-up rates for critical norms such as \eqref{e.quantbluprateTao} may enable to rule out certain blow-up scenarios for which the numerically computed growth of the $L^3$ norm is too slow. However,  the extremely slow triple logarithmic rate in \eqref{e.quantbluprateTao} means testing it may be beyond computing capacities, as is emphasized in Hou's recent paper \cite{hou2022potentially}.
\item As is well understood for dispersive equations, finding blow-up rates and concentration estimates is a first step toward understanding potential blow-up profiles.
\item In turbulence theory, cascade processes are the dominant feature of the inertial range, where nonlinear inertial effects dominate (or are in balance with) viscous dissipative effects. Below certain scales though, so-called `Kolmogorov scales', for very high wavenumbers or small spatial scales, dissipative effects dominate. Estimating those dissipative scales quantitatively is one of the main objectives of the quantitative regularity theory, see Section \ref{sec.quantkol} and Objective \ref{goal.quantdissip}. 
\end{enumerate}

\subsection*{Outline of the paper}

This survey paper is partly based on several talks given in the past two years. Further comments and topics are found in the habilitation thesis \cite{CP-hdr}. The first three sections are devoted to `weak concentration' (Section \ref{sec.weakconc}), `strong concentration' (Section \ref{sec.strongconc}) and a fundamental tool for the quantitative regularity, namely local-in-space smoothing (Section \ref{sec.locinspace}). The rest of the paper is concerned with quantitative regularity and blow-up rates for critical norms. Section \ref{sec.threefacets} describes three facets of quantitative regularity: blow-up rates, quantitative regularity estimates and quantitative estimates for dissipative scales. Moreover, two cases are studied, the case of $L^{5}_{t,x}$ which is a toy model, and the case of $L^\infty_tL^3_x$ studied by Tao. Section \ref{sec.genstrat} focuses on the strategy to estimate the dissipative scales in a quantitative manner. Section \ref{sec.quanttypeI} concentrates on the Type I case. Section \ref{sec.further} reviews some recent developments in the wake of \cite{Tao19}. Section \ref{sec.mildbreak} shows that the scaling barrier can be slightly broken thanks to good quantitative estimates in the critical case. This section is based on the paper \cite{BP21jmfm}. There a (partial) answer to a question asked in \cite{Tao19} about the blow-up of certain Orlicz norms is given. Finally in Section \ref{sec.summary}, we summarize some results in two tables.

Notice that $C$ and $c$ are universal positive constants that may change from line to line.

\section{Weak concentration}
\label{sec.weakconc}

In this paper we distinguish between:
\begin{description}
\item[`weak concentration'] of norms near a potential blow-up time $T^*$; these results assert the existence of points $x(t)$ (or a sequence of points $x_n$) and scales $\lambda(t)$ (or a sequence of scales $\lambda_n$) such that certain norms accumulate on $B_{x(t)}(\lambda(t))$ as $t\rightarrow T^*$;
\item[`strong concentration'] of norms near a potential space-time blow-up point $(x^*,T^*)$;\footnote{A `blow-up/singular point' $(x^*,T^*)$ is a point for which the solution is unbounded in any parabolic cylinder $Q_{(x^*,T^*)}(r)=B_{x^*}(r)\times(T^*-r^2,T^*)$ centered at that point. Conversely, a `regular point' is a point which is not a singular point. It is known that determining whether or not the singular points occur for Leray-Hopf solutions is equivalent to the Millenium problem.} these results assert the existence of scales $\lambda(t)$ (or a sequence of scales $\lambda_n$) such that certain norms accumulate on shrinking balls $B_{x^*}(\lambda(t))$ as $t\rightarrow T^*$.
\end{description}

\subsubsection*{\underline{A short review of concentration for certain nonlinear PDEs}}
The first results on concentration near potential singularities date back to more than 30 years ago. They mainly fall into the class of `weak concentration'. The study of `mass concentration' i.e. concentration of the $L^2$ norm for the nonlinear Schr\"odinger equations has triggered a lot of developments in this direction, in the wake of the seminal results by Weinstein \cite{weinstein1989nonlinear}, Merle and Tsutsumi \cite{MT90}, Nawa \cite[Theorem B and Theorem C]{Nawa90} and \cite[Theorem B]{Nawa92}, 
Merle \cite{Merle92,Merle93}. These results were followed by many others: Bourgain\cite{Bourgain98}, Nawa and Tsutsumi \cite{NT98}, Hmidi Keraani \cite[Corollary 1.8]{HK06} for nonlinear Schr\"odinger, Kenig, Ponce and Vega \cite[Corollary 1.4]{KPV00} for the KdV equation, Merle and Zaag \cite[Theorem 1, (ii)]{MZ05} for the semilinear wave equation\ldots\ We also refer to the books by Cazenave \cite[Section 6.5]{cazenave1989introduction}, Tao \cite{Tao06} and Sulem and Sulem \cite[Section 5.2.4 and Section 14.3.2]{sulem2007nonlinear}.

There are also a number of results that fall into the category of `strong concentration' results, especially for the nonlinear Schr\"odinger equation, for radially symmetric solutions: Merle and Tsutsumi \cite{MT90}, Tsutsumi \cite[Theorem 1.1]{Tsutsumi90}, Holmer Roudenko \cite[Theorem 1.2]{HR07}\ldots

The topic of concentration is also strongly tied to proving the blow-up of certain critical norms. The recent paper of Mizoguchi and Souplet on the semilinear heat equation states a strong concentration property for a Type I singularity \cite[Lemma 3.1]{MizoguchiSouplet19} that is key to proving the blow-up of a critical norm in that case; see also Miura and Takahashi \cite{MT22} without the Type I assumption.

As for fluids, we mention two results dealing with the possible energy concentration in solutions of Euler and Navier-Stokes equations. These results are a bit orthogonal to the concentration results on which we focus in this paper, but are also interesting directions. Chae and Wolf \cite{chae2020energy} proved that for 3D Euler under a Type I condition there is no concentration of the energy into isolated points at possible blow-up times. For 3D Navier-Stokes, Arnold and Craig \cite[Theorem 4.5]{AC10-concentration} gave a lower bound on the energy concentrating set.

\subsubsection*{\underline{Weak concentration for the Navier-Stokes equations}}

Li, Ozawa and Wang \cite[Theorem 1.2]{LOW18} pro\-ve what we believe is the first concentration result for potential blow-up solutions of the 3D Navier-Stokes equations. It is proved that for any Leray-Hopf solution $U$ that first blows-up at time $T^*\in(0,\infty)$, one has for any $q\in[1,\infty]$, the following concentration of the $L^q$ norm:
\begin{equation}\label{e.conclow}
\|U(\cdot,t)\|_{L^q(|\cdot-x_n|\lesssim\omega(t_n)^{-1})}\gtrsim\omega(t_n)^{1-\frac3q}
\end{equation}
where $t_n\rightarrow T^*$, $x_n\in\R^3$ and $\omega(t_n):=\|U(\cdot,t_n)\|_{L^\infty(\R^3)}\stackrel{t\rightarrow T^*}{\longrightarrow}\infty$. The proof is based on: (i) selecting a sequence of times $t_n$ on which one has large growth of the $L^\infty$ norm of the solution, (ii) showing that the low frequencies $\lesssim \omega(t_n)$ contribute to a large part of the growth of the $L^\infty$ norm of $U$. 

Let us make two comments on this result. First, it is a concentration result that holds for any Leray-Hopf solution that blows-up. The price to pay for this generality is the fact that the concentration is weak, i.e. not localized in space. Second, the concentration \eqref{e.conclow} holds for any norm in the Lebesgue scale, no matter whether the norm is subcritical ($q\in(3,\infty]$), critical ($q=3$) or supercritical ($q\in[1,3)$). Notice that in this last case, the lower bound in \eqref{e.conclow} goes to zero as $t\rightarrow T^*$. In the case of $q\in(3,\infty]$, one can bound the right hand side of \eqref{e.conclow} from below by $(T^*-t_n)^{-\frac12(1-\frac3q)}$ thanks to Leray's \cite{Leray} lower bound, see \eqref{e.lerayblup} below, and get the concentration on a ball of size $\sqrt{T^*-t_n}$.

Existence results of mild solutions with $L^q_{uloc}$ initial data also enable to prove weak concentration for potential blow-up solutions to the Navier-Stokes equations. This was remarked in \cite[Corollary 1.1]{MMP17a}, where the existence of mild solutions in $L^q_{uloc}$ is combined with a simple scaling argument to yield that for every $t\in (0, T^*)$, there exists $x(t)\in\mathbb{R}^3$, such that for any $q\in[3,\infty]$,\footnote{Let us stress that although the paper \cite{MMP17a} is actually concerned with the existence of mild solutions for data in $L^q_{uloc}$ in the half-space, we state \eqref{maekawamiuraprange} for the whole-space. This concentration result is implied by the existence result of mild solutions in the whole-space from \cite{maekawa2006} and the argument of \cite[Corollary 1.1]{MMP17a}.}
\begin{equation}\label{maekawamiuraprange}
\|U(\cdot, t)\|_{L^{q}(|\cdot-x(t)|\leq\sqrt{T^*-t})}\gtrsim (T^*-t)^{-\frac{1}{2}(1-\frac{3}{q})}.
\end{equation}
This strategy is robust enough to apply to the half-space $\R^3_+$ as in \cite{MMP17a}.

In \cite[Theorem 1.6 (i)]{KangMiuraTsai20-concL2}, Kang, Miura and Tsai prove a statement that can be read as weak concentration result for the supercritical $L^2$ norm near potential singularities of the Navier-Stokes equations: there exists $\gamma_{univ}\in(0,\infty)$, $S\in(0,\infty)$ and a function $x=x(t)\in\R^3$ such that for all $t\in(0,T^*)$,
\begin{equation}\label{ec7.concL2}
\frac1{\sqrt{T^*-t}}\int\limits_{B_{x(t)}\big(\sqrt{\frac{T^*-t}{S}}\big)}|U(x,t)|^2\, dx>\gamma_{univ}.
\end{equation}
This result is in the vein of the one of Bradshaw and Tsai \cite[Theorem 8.2]{bradshaw2022local}, see also Gruji{\'c} and Xu \cite[Theorem 4.1]{GrujicXu2019-dynrestr}.

The bottom line is that `weak concentration' results are in general a consequence of global regularity results or local well-posedness results and hence use a limited amount of specific structure of the equations. In order to have `strong concentration' results, one needs more localized regularity results at the price of additional a priori assumptions such as Type I. This is the topic of the next section.

\section{Local-in-space smoothing}
\label{sec.locinspace}

The idea behind local-in-space smoothing is very natural. Assume that one is given a rough initial data, for instance finite-energy, that happens to be more regular, in the sense critical or subcritical, on some fixed ball, say $B_0(1)$ to fix the ideas. 
For Navier-Stokes the general question becomes: 
\begin{quote}
Is the local smoothing due to the heat part of the equation strong enough to compensate for the nonlocal effects of the pressure that tend to propagate irregularities of the solutions from large-scale spatial scales to $B_0(1)$?
\end{quote}

The answer is yes in many situations, when the local `regular' data is taken in critical or subcritical Lebesgue, Lorentz or Besov spaces. More surprisingly, such results remain even true for the Navier-Stokes equations in the half-space with no-slip boundary condition, where nonlocal effects are known to be strong, see for instance \cite{Kang05,kang2021finite,SSv10}. When the answer is yes, we call such results `local-in-space short-time smoothing'. The solution $U$ then satisfies bounds of the type
\begin{equation}\label{e.bddLinftyU}
\sup_{t \in (0,S(M))} t^{\frac12} \|U(\cdot,t)\|_{L^\infty(B_0(\frac12)))} \lesssim 1
\end{equation}
under the condition that $\|U_0\|_{L^3(B_0(1))}\lesssim 1$ and for a short time $S=S(M)\ll 1$. 
This line of research was pioneered by Jia and \v{S}ver\'{a}k in the seminal paper \cite{JS14}. We find three main classes of methods that we sketch below. For more details and in particular statements of theorems, we refer to \cite[Theorem 3.1, Section 2 and 3]{JS14}, \cite[Theorem 1, Section 2 and 4]{BP18}, \cite[Chapter 5]{CP-hdr}, \cite[Theorem 1.1 and Section 3]{KangMiuraTsai20-concL2} and \cite[Theorem 1.6 and Section 5]{kwon2021role}.

Let us stress that local-in-space smoothing is a versatile tool that proved to be efficient in many situations, such as the existence proof of self-similar solutions \cite{JS14}, strong concentration (see below Section \ref{sec.strongconc}) and quantitative regularity (see below Section \ref{sec.threefacets} and Section \ref{sec.quanttypeI}). For the last point in particular, it is important to have quantitative versions of local-in-space smoothing, see for instance \cite[Theorem 5.1]{BP21cmp}.

Local-in-space smoothing also helps us gain understanding of a very natural question: for which initial data does the associated solution exhibit improved regularity properties?

\subsection{Lin type compactness methods}

We work with $U$ a global-in-time local energy solution to the Navier-Stokes equations on $\R^3\times (0,\infty)$ with initial data $U_0\in L^2_{uloc}(\R^3)$ with some mild decay at space infinity in order to rule out parasitic solutions and get a formula for the pressure.\footnote{The framework here is that of global solutions. However, it is possible to localize the results and state them for suitable solutions, see for instance \cite{KMT18,kwon2021role,ABP21}.} Assume that $U_0|_{B_0(1)}\in L^q(\R^3)$, with $q\in[3,\infty]$.\footnote{If the data is locally in the critical space $L^3$, one requires in addition smallness of the data.}

The compactness method proceeds in two steps. One first decomposes the solution $U$ to the Navier-Stokes equations into a mild solution $a$ originating from the critical or subcritical data $U_0|_{B_0(1)}$\footnote{One needs to properly extend this data as a compactly supported divergence-free function.} and a perturbation $V$ solving a perturbed Navier-Stokes system, with critical or subcritical drift terms and initial data locally zero in $B_0(1)$. The perturbation has a small energy locally in $B_0(\frac12)$ near the initial time and can be extended by zero backward in time. Hence the regularity of $V$ falls into the realm of epsilon-regularity results, which is the second step of this method.

The method for establishing the epsilon-regularity for the perturbed equation is inspired from the compactness method \cite{Lin} that Lin used to prove epsilon-regularity for the Navier-Stokes equations. Since the equation for the perturbation $V$ is a Navier-Stokes equation with drift terms, one needs to discriminate between subcritical and critical drifts:
\begin{itemize}
\item For subcritical drifts ($q\in(3,\infty]$), one has improved regularity for the limit equation in the compactness argument, hence one can directly prove local space-time H\"older regularity near initial time of the perturbation $V$. This was done in \cite{JS14}.
\item For critical drifts ($q=3$), there is no improved H\"older regularity at the limit in general. One therefore aims at first proving a subcritical Morrey bound for the perturbation that just misses boundedness; see \cite{KMT18}. Then one can combine this subcritical information for $V$ with subcritical information for the mild solution $a$ away from initial time to apply standard epsilon-regularity results for the Navier-Stokes equations to give boundedness of $V$ up to the initial time; see for instance \cite[Section 5.3]{KMT18}. It is also possible using information from the initial data to bootstrap the regularity of the perturbation to be H\"older continuous near the initial time; see for instance \cite[Section 3]{BP18}.
\end{itemize}

In the half-space with no-slip boundary condition, we were able to use the compactness method to prove local-in-space smoothing for subcritical and critical data in the Lebesgue scale; see \cite{ABP21}. This work relies on the new estimates for the harmonic pressure obtained in \cite{MMP17b}. The compactness scheme is convenient, because in the critical case, the smallness of the drift can be incorporated in the scheme, which avoids proving regularity for the limiting Stokes equation with drift terms, as is done for the whole-space in \cite[Lemma 2.2]{JS14} in the subcritical case or \cite[Section 4]{KMT18} in the critical case.

The bottom line is that the compactness method is flexible to handle both the global and the local settings and the regularity away or near boundaries.

\subsection{Caffarelli-Kohn-Nirenberg type methods}

The general scheme of the method is exactly as the one previous of the previous method: decomposition $U=a+V$ and smallness of the local energy of the perturbation $V$ (first step) and epsilon regularity for the perturbed Navier-Stokes equation with drift terms (second step). The difference is in the way the epsilon regularity is proved. In the paper \cite[Section 2]{BP18}, we prove an epsilon regularity result for the perturbed Navier-Stokes equations with critical drift terms by using an iteration scheme \emph{\`a la} Caffarelli, Kohn and Nirenberg \cite{CKN82}. The criticality of the drift $a$ associated to $L^3$ data causes difficulties, as was the case for the compactness method, and therefore only enables to propagate a subcritical Morrey bound for the perturbation $V$.

The main advantage of this method lies in the fact that given the current state of the art, it is the only method that manages to handle data locally small in the borderline endpoint Besov space $B^{-1+\frac3q}_{q,\infty}$ for $q\in(3,\infty)$; see \cite[Appendix C]{BP18}. The results can also be localized, see \cite[Theorem 3]{BP18}.

\subsection{Scaled local energy methods}

This approach was started by Kang, Miura and Tsai \cite{KangMiuraTsai20-concL2}. Roughly speaking, it relies on working with the scale-invariant energy and pressure
\begin{equation*}
E_r(t):=\frac1r\int\limits_{B_0(r)}|U(\cdot,t)|^2+\frac1r\int\limits_0^t\int\limits_{B_0(r)}|\nabla U|^2+\frac1{r^2}\int\limits_0^t\int\limits_{B_0(r)}|p|^\frac32,
\end{equation*}
and to propagate the smallness at initial time, i.e. of $\sup_{r\in(R,1)}E_r(0)$ where $R$ is a given non negative scale, forward-in-time via local energy estimates and a nonlinear Gronwall-type inequality for
$$
\mathscr E_{R,\hat R}(t)=\sup_{r\in[R,\hat R]}\sup_{s\in(0,t)}E_r(s),
$$
for a well-chosen parameter $\hat R$. 

There are two main advantages of this method. The first advantage is a technical one. The first step of the previous two methods, which consists in splitting the solution into $U=a+V$ where $a$ is the mild solution and $V$ the perturbation is not needed here any longer. One directly works with the solution. The second advantage is that it enables to directly prove local-in-space smoothing for small scaled local kinetic energy \cite[Theorem 1.1]{KangMiuraTsai20-concL2}.

All the local-in-space smoothing results in the local setting mentioned so far are for suitable solutions. Hence the pressure is a priori assumed to be in the Lebesgue space $L^\frac32$. Building upon the method of Kang, Miura and Tsai \cite{KangMiuraTsai20-concL2}, Kwon \cite[Theorem 1.6]{kwon2021role} was able to extend the local-in-space smoothing result for data locally in $L^3$ to the class of dissipative solutions, for which the pressure is barely a distribution. Moreover in \cite[Theorem 1.6]{kwon2021role} the time interval for which the local-in-space smoothing occurs is independent of the pressure.

\section{Strong concentration}
\label{sec.strongconc}

As mentioned above, see the preamble of Section \ref{sec.weakconc}, strong concentration is the accumulation of norms near potential blow-up points, on balls centered at a singularity, whereas `weak concentration' shows accumulation near blow-up times and is not well localized in space. 

Given the current state of the art, strong concentration can be proved in two cases:
\begin{itemize}
\item either under some symmetry assumption on the solution, such as radial symmetry for solutions of the nonlinear Schr\"odinger equation, see for instance \cite{MT90,Tsutsumi90,HR07},
\item or under a Type I assumption, 
see for instance \cite[Lemma 3.1]{MizoguchiSouplet19} for the semilinear heat equation.
\end{itemize}

In principle, strong concentration results for the Navier-Stokes equations are simply deduced from local-in-space smoothing results by a contradiction argument, on condition that one has a good workable notion of Type I.\footnote{The case of Type I a priori control contains many still unresolved blow-up Ans\"atze, such as Backward Discretely Self-Similar solutions (see \cite{CW17}).} We assume that the solution satisfies the following generalized Type I bound: for fixed $M,\ T^{*}\in (0,\infty)$ and a fixed radius $r_0\in (0,\infty]$,
\begin{align}
\begin{split}\label{ec6.gentypeI}
&\sup_{\bar x\in\R^3}\sup_{r\in(0,r_0)}\sup_{T^{*}-r^2<t<T^{*}}r^{-\frac12}\Bigg(\int\limits_{B_{\bar x}(r)}|U(x,t)|^2dx\Bigg)^\frac12\leq M.
\end{split}
\end{align}
In order to make sense of the condition \eqref{ec6.gentypeI} in the case when $r_0>\sqrt{T^*}$, we may extend $U$ by zero in negative times.\footnote{Notice that it is not immediatly clear that the generalized notion of Type I \eqref{ec6.gentypeI} is implied by more classical notions of Type I, such as the ODE blow-up Type I
\begin{equation}\label{e.typeIa}
\sqrt{T^{*}-t}|U(x,t)|\leq M'
\end{equation}
It turns out to be true, see \cite{SerZajac} and the review article \cite[pages 844-849]{Seregin2018}. This makes \eqref{ec6.gentypeI} a good notion. In the half-space with no-slip boundary conditions such an implication remains true, but is much harder to prove due to the strong nonlocality of the pressure. It is proved in \cite{BP20geom}.} That condition ensures that blow-up profiles belong to the class of local energy solutions in which local-in-space smoothing results are proved, see Section \ref{sec.locinspace}. 

\begin{figure}[t]
\begin{center}
\includegraphics[scale=.55]{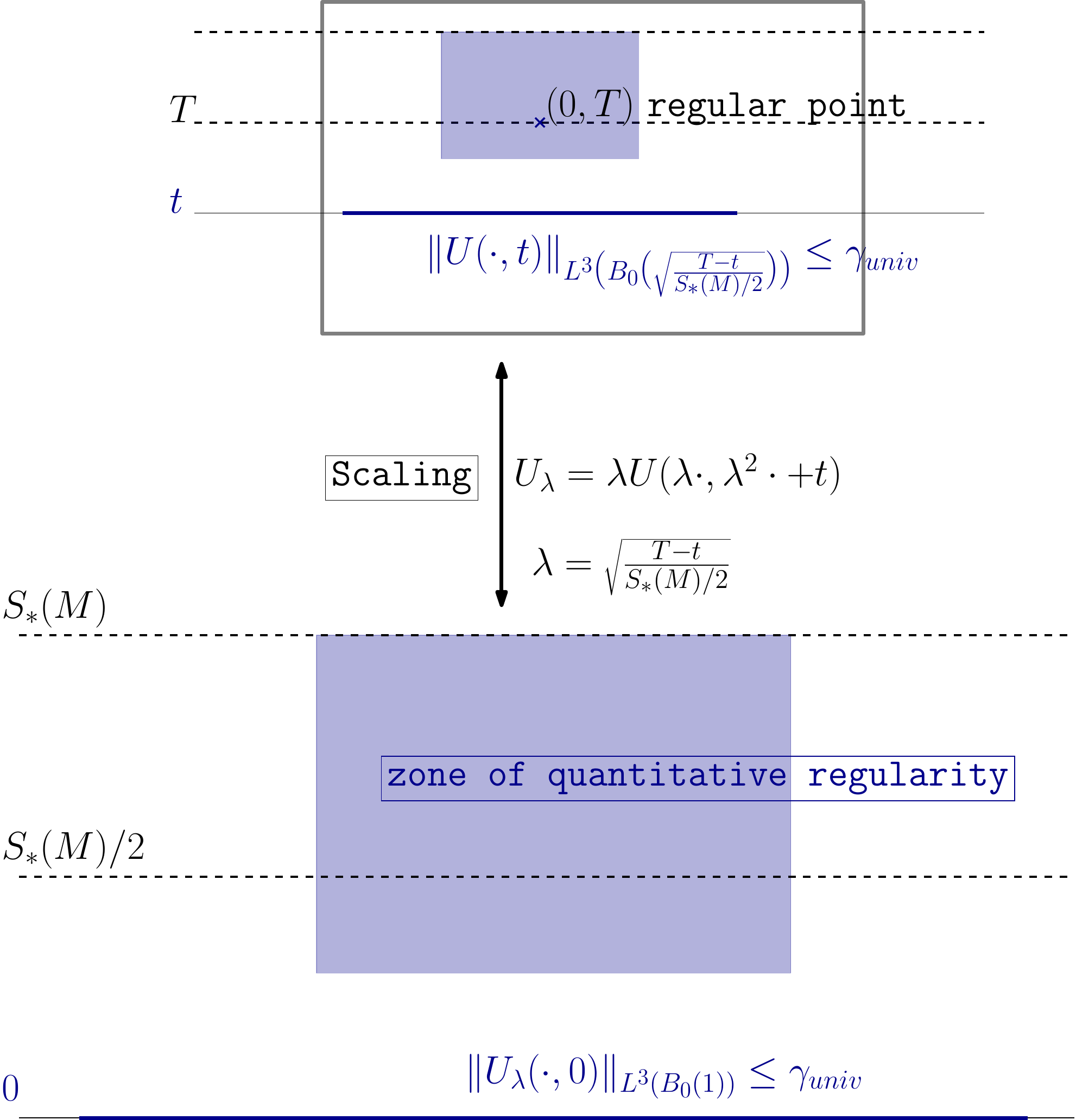}
\caption{Quantitative local-in-space short-time smoothing and concentration near potential Type I singularities}
\label{fig.c6locsmoothing}
\end{center}
\end{figure}

In \cite{BP18}, strong concentration of the critical $L^3$, $L^{3,\infty}$ and Besov norms $B^{-1+\frac3q}_{q,\infty}$ ($q\in(3,\infty)$) are proven. Figure \ref{fig.c6locsmoothing} explains the idea to deduce strong concentration of the $L^3$ norm from local-in-space smoothing for critical $L^3$ data. For a solution $U$ satisfying the Type I bound \eqref{ec6.gentypeI}, first blowing-up at $T^*$ and such that the space-time point $(0,T^*)$ is a singularity, we have
\begin{equation}\label{e.conblup}
\|U(\cdot,t)\|_{L^3\big(B_0\big(\sqrt{\frac{T^{*}-t}{S(M)/2}}\big)\big)}>\gamma_{univ},
\end{equation}
for all $t\in (t_*(T^*,M,r_0),T^{*})$ and $S(M)$ a time appearing in the local-in-space smoothing estimate \eqref{e.bddLinftyU}. 

In \cite[Theorem 1.6 (ii)]{KangMiuraTsai20-concL2} Kang, Miura and Tsai prove the strong concentration of the scaled energy on concentrating balls. Their result holds under the generalized Type I condition \eqref{ec6.gentypeI}. There exists $\gamma_{univ}\in(0,\infty)$ and $S(M)\in(0,\infty)$ such that for all $t\in(t_*(T^*,M,r_0),T^*)$,
\begin{equation}\label{ec7.concL2typeI}
\frac1{\sqrt{T^*-t}}\int\limits_{B_{0}\big(\sqrt{\frac{T^*-t}{S(M)}}\big)}|U(x,t)|^2\, dx\gtrsim 1.
\end{equation}
Notice that under the stronger ODE blow-up Type I condition \eqref{e.typeIa}, this result simply follows from \eqref{e.conblup} and interpolation.

Figure \ref{fig.graphconc} on page \pageref{fig.graphconc} summarizes some results about weak and strong concentration for the Navier-Stokes equations, as well as local-in-space smoothing.

The bottom line is that strong concentration of certain scale-invariant quantities is at the heart of the quantitative regularity, in particular for proving quantitative estimates of dissipative scales; see Section \ref{sec.quantkol}, Section \ref{sec.genstrat} and Section \ref{sec.typeIestdissscales} below. From now on the paper is concerned with such questions relating to quantitative regularity.

\section{Three facets of quantitative regularity}
\label{sec.threefacets}

In our view, there are three main objectives: (i) quantitative blow-up rates, see Objective \ref{obj.bluprate} below, (ii) quantitative regularity estimates, see Objective \ref{obj.quantreg} below, (iii) quantitative estimates of dissipative scales, see Objective \ref{goal.quantdissip} below. The last objective is in some sense more fundamental, since the two others will in general follow from it. We formulate the general objectives in a somewhat loose way. In the rest of the paper, we will explain how the examples fit into that abstract framework. In particular in this section we illustrate the objectives in two cases: the case of a priori control in $L^5(\R^3\times(-1,0))$ which is a critical non-borderline space,\footnote{In this case the qualitative regularity is the classical Lady\v{z}enskaja-Prodi-Serrin criteria, see for instance \cite{Serrin1962,T90}.} and the case of a priori control in $L^\infty(-1,0;L^3(\R^3))$\footnote{In this case the qualitative regularity is the result of Escauriaza, Seregin and \v{S}ver\'{a}k \cite{ESS2003}.} which is a critical borderline space.

\subsection{Quantitative blow-up rates}

It is known since the seminal work of Leray \cite{Leray} that subcritical norms blow-up with an algebraic rate near a potential blow-up time $T^*$, i.e.
\begin{equation}\label{e.lerayblup}
\|U(\cdot,t)\|_{L^q(\R^3)}\gtrsim (T^*-t)^{-\frac32(\frac13-\frac1q)}\quad\mbox{for}\ q\in(3,\infty)\ \mbox{and}\ t\in[0,T^*).
\end{equation}

\begin{goala}\label{obj.bluprate}
Show that for certain critical spaces $L^p_t\mathcal A_x$ and $\mathcal B$,\footnote{Here $\mathcal A,\, \mathcal B$ are certain Banach spaces contained in $\mathcal S'(\R^3)$ and $p\in[1,\infty]$. Notice that $\mathcal B$ is scaling-invariant, as well as $L^p_t\mathcal A_x$.} 
there exists an explicit positive function $\mathscr F$ or $\mathscr F_{p,\mathcal A}$ on $[0,\infty)$ such that $\mathscr F(s),\ \mathscr F_{p,\mathcal A}(s)\stackrel{s\rightarrow 0}\longrightarrow \infty$ and for a smooth solution with enough decay\footnote{Under this assumption the solution is smooth on $(0,T)$ for any $T<T^*$, hence is a classical solution; see \cite[Section 1.4]{BP21cmp} for a definition.} $U$ to the Navier-Stokes equations on $\R^3\times(0,T^*)$ blowing-up at time $T^*$,
\begin{equation}\label{e.quantbluprate}
\frac{\|U\|_{L^p(0,t;\mathcal A)}}{\mathscr F(T^*-t)}\gtrsim 1\quad\mbox{when}\ t\ \mbox{is sufficiently close to}\ T^*.
\end{equation}
or 
\begin{equation}\label{e.quantblupratebis}
\frac{\|U(\cdot,t)\|_{\mathcal B}}{\mathscr F_{p,\mathcal A}(T^*-t)}\gtrsim 1\quad\mbox{when}\ t\ \mbox{is sufficiently close to}\ T^*.
\end{equation}
\end{goala}

\subsubsection*{\underline{The toy model $L^5_{t,x}$}}

Using the same reasoning as \eqref{e.quanttoblup}, the quantitative estimate \eqref{L5quantest} combined with Leray's blow-up rate \eqref{e.lerayblup} imply that
\begin{equation*}
\|U\|_{L^5(0,t;L^5(\R^3))}\gtrsim \sqrt{T^*}(\log(T^*-t))^\frac15,\quad\mbox{for all}\quad t\in(\tfrac{T^*}2,T^*).
\end{equation*} 
Not surprisingly, this estimate is compatible with (and actually a consequence of) the Leray blow-up rate \eqref{e.lerayblup} in the case $q=5$.

\subsubsection*{\underline{The case of $L^\infty_tL^3_{x}$}}

The critical borderline case $q=3$ was open until 2019 and a remarkable paper of Tao \cite{Tao19}. There Tao proves the rate \eqref{e.quantbluprateTao}. Notice that this rate follows from the lower bound \eqref{e.quantbluprate} proved in \cite{Tao19} with $p:=\infty$, $\mathcal A:=L^3(\R^3)$ and $\mathscr F(s):=(\log\log\log s)^{c}$, 
for a universal constant $c\in(0,\infty)$.

\subsection{Quantitative regularity estimates}

\begin{goala}\label{obj.quantreg}
Show that for certain critical spaces $L^p(-1,0;\mathcal A)$ and $\mathcal B$,\footnote{\label{foot.critAB}Here $\mathcal A,\, \mathcal B$ are certain Banach spaces contained in $\mathcal S'(\R^3)$ and $p\in[1,\infty]$. Notice that $\mathcal B$ is scaling-invariant, as well as $L^p(-1,0;\mathcal A)$.}
there exists an explicit positive function $\mathscr G$ on $[0,\infty)$,\footnote{It is hard to give a general statement covering all results in this line of research. Notice that Theorem B in \cite{BP21cmp} is a quantitative regularity result without a priori control of a global scale-invariant quantity. Such a result takes the following form
$$
\|U\|_{L^{\infty}(\mathbb{R}^3\times (-t_*,0))}\leq \mathscr G\big(\sup_{t_j\nearrow 0}\|U(\cdot,t_j)\|_{L^3(\R^3)}\big),
$$
where $t_*=t_*(\sup_{t_j\nearrow 0}\|U(\cdot,t_j)\|_{L^3(\R^3)})$. 
That result quantifies Seregin's 2012 liminf qualitative criteria \cite{seregin2012}, and is hence a result that goes beyond the critical case. For more details about \cite[Theorem B]{BP21cmp}, we also refer to \cite[Theorem 6.3 and Section 6.3.4]{CP-hdr}. Furthermore, in the same vein, we get regularity at time $t=0$ for large $L^3$ data at $t=-1$ if the profile at time $t=0$ is quantitatively small; see \cite[Theorem 4.1 (i)]{AlbrittonBarkerBesov2018} for a qualitative statement and \cite[Proposition 4.4]{BP21cmp} for a quantitative statement.} such that for a critically bounded smooth solution $U$ (with enough decay) of the Navier-Stokes equations on $\R^3\times(-1,0)$,
\begin{equation}\label{quantgenform}
\|U\|_{L^{\infty}(\mathbb{R}^3\times (-\frac{1}{2},0))}\leq \mathscr G(\|U\|_{L^p(-1,0;\mathcal A)},\|U(\cdot,0)\|_{\mathcal B}).
\end{equation} 
\end{goala}

We previously remarked that certain qualitative regularity results in terms of critical norms can be quantified abstractly, see \cite[Introduction]{BP21cmp}. This is the case of the Escauriaza, Seregin and \v{S}ver\'{a}k \cite{ESS2003} result that can be abstractly quantified via the use of the persistence of singularities \cite{rusin2011minimal}. The focus of this survey paper is to derive an explicit formula for $\mathscr G$. 

Let us also note that a general quantitative result such as \eqref{quantgenform} can be combined with the Leray blow-up rate \eqref{e.lerayblup} to yield quantitative blow-up rates of the form \eqref{e.quantbluprate}. Indeed, by scaling\footnote{We recall that the Navier-Stokes equations are invariant under the scaling $U_\lambda=\lambda U(\lambda\cdot,\lambda^2\cdot)$, for $\lambda>0$. Here we take $\lambda=\sqrt{t}$.}
\begin{align}\label{e.quanttoblup}
\begin{split}
\frac1{\sqrt{t}}\mathscr G(\|U\|_{L^\infty(0,t;L^3(\R^3))})&\ =\frac1{\sqrt{t}}\mathscr G(\|U_{\sqrt{t}}\|_{L^\infty(0,1;L^3(\R^3))})\\
&\ \geq\frac1{\sqrt{t}}\|U_{\sqrt{t}}\|_{L^\infty(\R^3\times(\frac 12,1))}\\
&\ =\|U\|_{L^\infty(\R^3\times(\frac t2,t))}\\
&\ \geq C(T^*-t)^{-\frac12}.
\end{split}
\end{align}

\subsubsection*{\underline{The toy model $L^5_{t,x}$}} 

Using the strategy described in Section \ref{sec.genstrat} below,\footnote{A strategy based on energy estimates on the level of the vorticity equation is described in \cite[Section 1.2.2]{BP21cmp}. It also yields a single exponential bound in $L^5(\R^3\times(-1,0))$, but has an additional dependence in terms of the initial data for the vorticity $\omega(\cdot,-1)$. Notice that $L^5$ energy estimates can be applied directly, see \cite{MontgomerySmith05-mildbreak}.} one can prove the following quantitative estimate:
\begin{equation}\label{L5quantest}
\|U\|_{L^{\infty}(\mathbb{R}^3\times (-\frac12,0))}\lesssim \exp\big(\|U\|_{L^5(\R^3\times (-1,0))}^5\big),
\end{equation}
where $C\in(0,\infty)$ is a universal constant. Here 
$$p=5,\quad\mathcal A=L^5(\R^3)\quad\mbox{and}\quad\mathscr G(A,B)\simeq\exp(A^5)$$ in the notation of the general form estimate \eqref{quantgenform}.

\subsubsection*{\underline{The case of $L^\infty_tL^3_{x}$}}

Tao shows that for classical 
solutions to the Navier-Stokes equations on 
belonging to the critical space $L^{\infty}(-1,0; L^{3}(\mathbb{R}^3))$,
\begin{equation}\label{taomainest}
\|U\|_{L^{\infty}(\mathbb{R}^3\times (-\frac12,0))}\lesssim \exp\exp\exp\big(\|U\|_{L^{\infty}(-1,0;L^3(\mathbb{R}^3))}^{c}\big),
\end{equation}
where $c\in(0,\infty)$ is a universal constant. Here 
$$p:=\infty,\quad\mathcal A:=L^3(\R^3)\quad\mbox{and}\quad\mathscr G(A,B)\simeq\exp(\exp(\exp(A^c)))$$ in the notation of the general form estimate \eqref{quantgenform}.

\subsection{Quantitative estimates of dissipative/Kolmogorov scales}
\label{sec.quantkol}

\begin{goala}\label{goal.quantdissip}
Show that for certain critical spaces $L^p(-1,0;\mathcal A)$ and $\mathcal B$,\footnote{As in footnote \ref{foot.critAB}, $\mathcal A,\, \mathcal B$ are certain Banach spaces contained in $\mathcal S'(\R^3)$ and $p\in[1,\infty]$.} there exists an explicit positive function $\mathscr H$ on $[0,\infty)$ such that for a critically bounded smooth solution with enough decay $U$ to the Navier-Stokes equations on $\R^3\times (-1,0)$, 
dissipative effects take over the nonlinearity\footnote{This is on a formal level. On a practical level, this is where well-posedness theory or epsilon regularity takes over and gives regularity.} for physical scales $$\lambda\leq \mathscr H(\|U\|_{L^\infty(-1,0;\mathcal A)},\|U(\cdot,0)\|_{\mathcal B})$$ or Fourier scales $$N\geq \mathscr H
(\|U\|_{L^\infty(-1,0;\mathcal A)},\|U(\cdot,0)\|_{\mathcal B}).$$
\end{goala}

The threshold between the scales where diffusive effects dominate vs. where nonlinear effects dominate is measured by the smallness vs. concentration of certain scale-critical quantities $\mathscr  S$. Once this threshold is estimated in a quantitative way, regularity criteria (epsilon-regularity, local smoothing\ldots) imply the quantitative bounds stated in Objective \ref{obj.bluprate}. The estimate of dissipative scales is hence the heart of the matter.

Figure \ref{fig.summaryquantL5} summarizes the general strategy in the three main cases described in this survey paper.

\begin{figure}[h]
\begin{center}
\begin{tabular}{|C{5cm}|C{6cm}|C{5cm}|}
\hline
\multicolumn{3}{|c|}{the global scale-critical standing assumptions}\\
$U\in L^5_{t,x}$ & $U\in L^\infty_tL^{3,\infty}_x$ and $U(\cdot,0)\in L^3$ & $U\in L^\infty_tL^3_x$\\
see Section \ref{sec.quantkol} & see Section \ref{sec.typeIestdissscales} & see Section \ref{sec.quantkol}\\
\multicolumn{3}{|c|}{prevent the following scale-critical quantities}\\
${\displaystyle(-t)^\frac15\|U(\cdot,t)\|_{L^5}}$ & ${\displaystyle\sqrt{-t}\int\limits_{B_{O\big(\|U\|_{L^\infty_tL^{3,\infty}_x}\sqrt{-t}\big)}}|\omega(\cdot,t)|^2}$  & ${\displaystyle N^{-1}|P_NU(t,x)|}$\\
defined by \eqref{e.estdissipL5} & defined by \eqref{vortconcintro} & defined by \eqref{e.freqbubbles}\\
\multicolumn{2}{|c|}{to concentrate to close to final time} & to concentrate for large frequencies \\
\multicolumn{2}{|c|}{i.e. for times $t_*<t<0$} & $N\geq N_*$\\
\hline
\multicolumn{3}{|c|}{no concentration, i.e. smallness, implies regularity}\\
  \hline
  \end{tabular}
\end{center}
\caption[Quantitative estimates via concentration]{Quantitative regularity via concentration: a summary}
\label{fig.summaryquantL5}
\end{figure}

\subsubsection*{\underline{The toy model $L^5_{t,x}$}} For solutions of the Navier-Stokes equations critically bounded in $L^5(\R^3\times(-1,0))$, we work with the following scale-invariant Weissler-Kato type norms
\begin{equation}\label{e.weisslerkato}
(-t)^\frac15\|U(\cdot,t)\|_{L^5(\R^3)},\quad\mbox{for}\quad t\in(-1,0).
\end{equation}
If the quantity defined by \eqref{e.weisslerkato} is small  
for all sufficiently small times $0>t>t_*(A)$, then Caffarelli, Kohn and Nirenberg type epsilon-regularity results imply the quantitative bound \eqref{L5quantest}.\footnote{Notice that well-posedness theory could also be used in replacement of epsilon-regularity.}

\begin{description}
\item[Objective \ref{goal.quantdissip} in the $L^5_{t,x}$ case] 
If the following statement\footnote{The $\ep$ in \eqref{e.estdissipL5} comes from epsilon-regularity criteria.}
\begin{align}\label{e.estdissipL5}
\mathscr S(t):=(-t)^\frac15\|U(\cdot,t)\|_{L^5(\R^3)}<\ep
\end{align}
fails for a certain $t\in(-1,0)$, find a quantitative upper bound $t_*\in(-1,0)$ for $t$.
\end{description}

In can be shown, applying the general strategy outlined in Section \ref{sec.genstrat}, that 
\begin{equation}\label{e.estt*L5}
-t>\frac12\exp\big(-\tfrac{32\|U\|_{L^5(\R^3\times(-1,0))}^5}{\ep^5}\big)=:-t_*(A).
\end{equation}
Hence, $\mathscr H(A,B):=-t_*(A)=\exp(-\frac{32A^5}{\ep^5})$ with the notation of Objective \ref{goal.quantdissip}.

\subsubsection*{\underline{The case of $L^\infty_tL^3_{x}$}}

Tao works with globally defined quantities due to a Fourier based approach. His analysis is based on the following scale-invariant quantities
\begin{equation}\label{e.freqbubbles}
\mathscr S(N;x,t):=N^{-1}|P_{N}U(x,t)|,\quad\mbox{for}\quad (x,t)\in\R^3\times(-1,0),
\end{equation}
where $P_N$ is the Littlewood-Paley projection on the frequency $N\in(0,\infty)$. Indeed, if the quantity $\mathscr S(N;x,t)$ defined by \eqref{e.freqbubbles} is small in terms of $A:=\|U\|_{L^{\infty}_{t}L^{3}_{x}(\mathbb{R}^3\times (-1,0))}$, uniformly in $(x,t)\in\R^3\times (-\frac12,0)$ and for high frequencies $N\geq N_*(A)$, then $\|U\|_{L^{\infty}_{x,t}(\mathbb{R}^3\times (-\frac{1}{8},0))}$ can be estimated explicitly in terms of $A$ and $N_{*}$. Related observations were made previously by Cheskidov and Shvydkoy \cite{cheskidov2010regularity,CS14} and Cheskidov and Dai \cite{ChesDai19}, but without the bounds explicitly stated. There the frequency $N_*$ is called the Kolmogorov scale and denoted $\Lambda$. If $N^{-1}\|P_{N}U\|_{L^{\infty}_{x,t}}\ll 1$ is small, then $N\|P_NU\|^2_{L^\infty_{t,x}}\ll N^2\|P_NU\|_{L^\infty_{t,x}}$ so that the diffusion dominates the nonlinearity, heuristically at least, since some of the frequencies in the paraproduct are neglected, see \cite{tao2013localisation}. 

From this perspective, Tao's aim is the following.
\begin{description}
\item[Tao's Objective \ref{goal.quantdissip}] Assume $$A:=\|U\|_{L^{\infty}_{t}L^{3}_{x}(\mathbb{R}^3\times (-1,0))}<\infty.$$ 
If the following statement
\begin{equation}\label{scaleinvarsmallhighfreqbis}
N^{-1}\|P_{N}U\|_{L^{\infty}_{x,t}(\mathbb{R}^3\times (-\frac{1}{2},0))}<\varepsilon(A)\,\,\,\textrm{for}\,\,\textrm{all}\,\,\,N\geq N_{*},
\end{equation}
fails for $\varepsilon(A)= A^{-c}$ and a certain frequency $N$, find a quantitative upper bound $N_*(A)$ for $N$.
\end{description}

In Tao's paper \cite[Theorem 5.1]{Tao19}, it is shown that 
\begin{equation}\label{e.estTaoN*}
\mathscr H(A,B):=N_*(A)\simeq\exp\exp\exp(A^{c})
\end{equation} 
with the notations of Objective \ref{goal.quantdissip}.

\section{A general strategy for quantitative estimates of dissipative scales}
\label{sec.genstrat}

In this part we describe a general strategy to fulfill Objective \ref{goal.quantdissip}, i.e. to estimate quantitatively the dissipative scales. The scheme is based on the study of the concentration of certain scale-invariant quantities $\mathscr S$, such as the Weissler-Kato norm \eqref{e.estdissipL5} in the case of $L^5_{t,x}$ or the frequency bubbles \eqref{e.freqbubbles} in the case of $L^\infty_tL^3_x$. One can summarize the idea as follows: 
\begin{enumerate}[label=(Step-\arabic*)]
\item \textbf{Propagation of concentration}\\
If a certain scale-invariant quantity $\mathscr S$ concentrates at a given scale, then it will concentrate at many disjoint scales. One relies on various tools such as bilinear estimates or local-in-space smoothing when applying quantitative Carleman inequalities. For this, one needs to work in space-time regions where one has good quantitative regularity pro\-perties: epochs and annuli of quantitative regularity.
\item \textbf{Summation of scales and coercivity of the standing assumption}\\
The standing critical assumption 
$$
\|U\|_{L^\infty(-1,0;\mathcal A)}<\infty\quad\mbox{or}\quad\|U(\cdot,0)\|_{\mathcal B}<\infty
$$
implies an upper bound on the number of disjoint scales where the scale-invariant quantity $\mathscr S$ concentrates, which directly yields a quantitative estimate of the estimate of dissipative scales $\mathcal H$ in Objective \ref{goal.quantdissip}.
\end{enumerate}

\subsection{The toy model $L^5_{t,x}$}

Let us now sketch the proof in the toy model case. Assume that $$A:=\|U\|_{L^5(\R^3\times (-1,0))}<\infty$$
and that \eqref{e.estdissipL5} fails for a certain $t$.  
In order to estimate the dissipative scales quantitatively, see Objective \ref{goal.quantdissip} above, we argue in the following two-step way:
\begin{enumerate}[label=(Step-\arabic*)]
\item \textbf{Backward propagation of Weissler-Kato norm concentration}\\
Using bilinear estimates for the Oseen tensor \cite{cannone1997generalization}  we get
\begin{equation*}
\|U(\cdot,t)\|_{L^5(\R^3)}>\frac{\ep}{(-t)^{\frac15}}\ \Rightarrow\ \|U(\cdot,t')\|_{L^5(\R^3)}>\frac{\ep/2}{(-t')^{\frac15}},\quad \mbox{for all}\quad t'\in(-1,2t).
\end{equation*}
\item \textbf{Summation of scales and coercivity of the standing assumption}\\
Integrating the concentration for $t'\in(-1,2t)$,
\begin{equation*}
\|U\|_{L^5(\R^3\times(-1,0))}^5\geq 
\int\limits_{-1}^{2t}\|U(\cdot,t')\|_{L^5(\R^3)}^5\, dt'
\geq -\frac{\ep^5}{32}\log(-2t).
\end{equation*}
Hence, we get the estimate \eqref{e.estt*L5}, which is the upper bound for the times where the Weissler-Kato norm concentrates.
\end{enumerate}

\begin{figure}[h]
\begin{center}
\includegraphics[scale=.55]{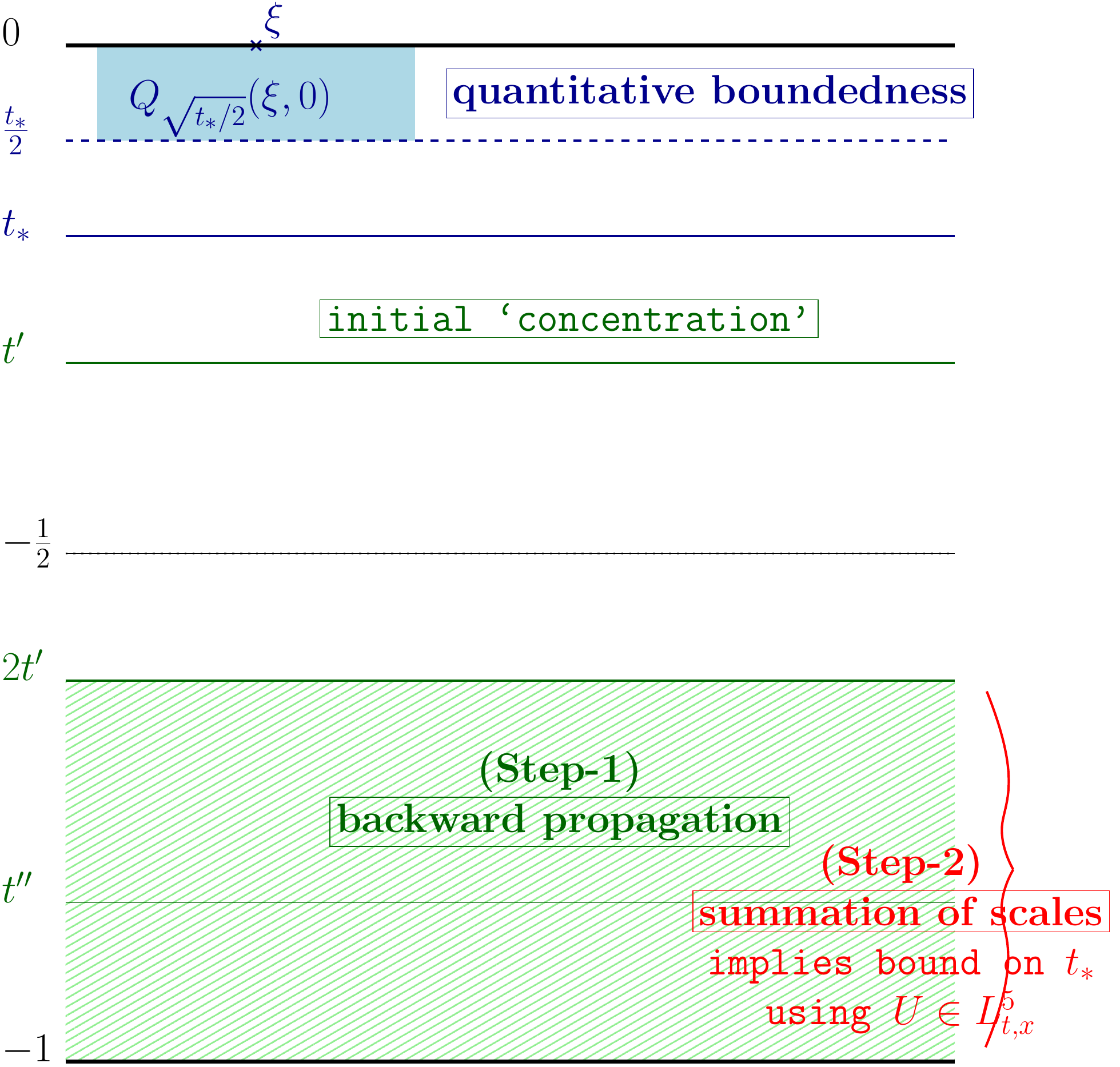}
\caption[Quantitative regularity via concentration]{Quantitative regularity via concentration of scale-critical quantities: the toy model $L^5_{t,x}$}
\label{fig.c7toy}
\end{center}
\end{figure}

\subsection{The case of $L^\infty_tL^3_x$}
\label{sec.taosstrat}

We now sketch the strategy of Tao to prove his main Objective \ref{goal.quantdissip} described on page \pageref{scaleinvarsmallhighfreqbis}. A slightly different summary of the strategy is given in \cite[Section 1.1]{BP21cmp} and in \cite[Section 1]{Tao19}. Assume that 
\begin{equation}\label{e.standingTao}
A:=\|U\|_{L^{\infty}(-1,0;L^{3}_{x}(\mathbb{R}^3))}<\infty
\end{equation}
and that \eqref{scaleinvarsmallhighfreqbis} fails for $\varepsilon(A)= A^{-c}$ and a certain frequency $N$. Hence there exists $(x_0,t_0)\in\R^3\times(-\frac12,0)$ such that 
\begin{equation}\label{frequencyconc}
N^{-1}|P_{N}u(x_0,t_0)|> A^{-c}.
\end{equation}
We call this the `initial concentration'.
\begin{enumerate}[label=(Step-\arabic*)]
\item\label{step.step1Tao} \textbf{Propagation of concentration}\\
The idea is to transfer the initial concentration \eqref{frequencyconc} in space and time in order to get a lower bound on the $L^3(\R^3)$ norm at time zero. That propagation of concentration, or `frequency bubbling', relies on:
\begin{enumerate}[label=(\roman*)]
\item \emph{Backward frequency bubbling}\\
For all $n\in\mathbb{N}$, there exists  a frequency $N_{n}\in(0,\infty)$, $(x_{n}, t_{n})\in\mathbb{R}^3\times (-1,t_{n-1})$ such that 
\begin{equation}\label{frequencyconciterate}
N_{n}^{-1}|P_{N_{n}}u(x_{n},t_{n})|> A^{-c}
\end{equation}
with 
\begin{equation}\label{parabolicdependence}
x_{n}=x_{0}+O((-t_{n})^{\frac12}),\quad N_{n}\simeq |-t_{n}|^{-\frac{1}{2}}.
\end{equation}
\item \emph{Transfer of concentration in Fourier space to physical space}\\ 
In order to use quantitative Carleman inequalities, see the next point, one needs to transfer the information on the concentration in Fourier space to physical space quantities, namely a scale-invariant enstrophy. To do this it seems important to have a priori control a global scale-invariant norm, such as \eqref{e.standingTao}.
\item\label{step.step1Taoii} \emph{Large-scale and forward-in-time propagation of concentration}\\
Using quantitative versions of the Carleman inequalities in \cite{ESS2003}, see  \cite[Proposition 4.2 and Proposition 4.3]{Tao19}, Tao shows that the lower bounds on the scale-invariant enstrophy can be transferred to a lower bound on the $L^{3}$ norm of $U$ at the final moment of time $0$. The applicability of the quantitative Carleman inequalities to the vorticity equation requires the `epochs of regularity' in the previous step and the existence of `good spatial annuli' where the solution enjoys good quantitative estimates. Specifically, Tao shows that for certain admissible time scales $S$,\footnote{These admissible time scales are related to where one has backward concentration in frequency space, see step `Backward frequency bubbling' above and \cite[Proof of Theorem 5.1]{Tao19}.} one has the concentration of the $L^3$ norm on the annulus $\{S^\frac12\leq|\cdot-x_{0}|\leq \exp(A^c)S^\frac12\}$, i.e.  \begin{equation}\label{L3lowerTao}
\int\limits_{S^\frac12\leq|x-x_{0}|\leq \exp(A^c)S^\frac12} |U(x,0)|^3 dx\geq \exp(-\exp(A^{c})).
\end{equation}
\end{enumerate}
\item\label{step.step2Tao} \textbf{Summation of scales and coercivity of the standing assumption}\\
Summing \eqref{L3lowerTao} over admissible time scales $S$ such that the annuli $\{S^\frac12\leq|\cdot-x_{0}|\leq \exp(A^c)S^\frac12\}$ are disjoint, one eventually obtains
\begin{equation*}
A^3\geq\int\limits_{\R^3}|U(x,0)|^3dx\gtrsim\log(N)\exp(-\exp(A^{c})).
\end{equation*}
This concludes the proof of the triple exponential upper bound \eqref{e.estTaoN*} for $N$.\footnote{Notice that due to the stacking of scales, we always get at least one exponential in the quantitative estimates.}
\end{enumerate}

Let us note that ideas in a similar spirit were used by Merle and Rapha\"el in \cite{MR06} to prove a quantitative rate of blow-up of a critical norm for radial solutions of a supercritical nonlinear Schr\"odinger equation; see \cite[Theorem 2]{MR06}. In particular a lower bound on annuli analogous to \eqref{L3lowerTao} is obtained and the final $\log$ blow-up rate is obtained by a similar stacking of scales argument; see \cite[Section 4.2]{MR06}.

\section{Quantitative regularity in the Type I case}
\label{sec.quanttypeI}

We focus here on the quantitative regularity results in the Type I case proved in \cite{BP21cmp}. Let us emphasize that the regularity in the Type I case is proved only for axisymmetric flows\footnote{In the half-space, this appears to still be an open problem even without swirl.} \cite{chen2008lower,chen2009lower,koch2009liouville}. We now show how the results in \cite{BP21cmp} fulfill the three objectives stated in Section \ref{sec.threefacets}.

\subsection{Three facets of the quantitative regularity in the Type I case}

\subsubsection*{\underline{Answer to Objective \ref{obj.bluprate}}}

\begin{resulta}[{\cite[Theorem A]{BP21cmp}}, localized quantitative rate of blow-up]\label{res.bluprateTypeI}
Assume that $U$ is a 
smooth solution with sufficient decay to the Navier-Stokes equations and that $T^*$ is a first blow-up time.\\
Assume in addition that $(0,T^*)$ is a Type I singular point i.e.
\begin{equation*}
A:=\|U\|_{L^{\infty}(0,T^*;L^{3,\infty}(\mathbb{R}^3))}<\infty.
\end{equation*}
Then the above assumptions imply that there exists  $S(A)\simeq A^{-30}$ such that for any $t\in (\frac{T^*}{2},T^*)$ 
and
\begin{equation}\label{theoARconstraint}
R\in \Big(\sqrt{\tfrac{T^*-t}{S(A)}},e^{A^{1022}}\sqrt{T^*}\Big)
\end{equation} 
we have
\begin{equation}\label{L3localisedlog}
\int\limits_{|x|<R} |U(x,t)|^3 dx\geq
\frac{\log\big(\frac{R^2}{A^{802}(T^*-t)}\big)}{\exp(\exp(A^{1025}))}.
\end{equation}
\end{resulta}

Estimate \eqref{L3localisedlog} is written in a scale-invariant form. Notice that this estimate implies an estimate of the form \eqref{e.quantbluprate}. Indeed taking $R=(T^*-t)^{\frac12-\delta}$ for $\delta>0$ and small, we have \eqref{e.quantbluprate} with $\mathcal B:=L^3(B_0((T^*-t)^{\frac12-\delta})))$, $p:=\infty$, $\mathcal A:=L^{3,\infty}(\R^3)$, $\mathscr F_{p,\mathcal A}\simeq-\log\big(A^{802}(T^*-t)^{2\delta}\big)$ and $|T^*-t|\lesssim_{A}1$.

Notice that contrary to \eqref{e.quantbluprateTao}, it is stated in \eqref{L3localisedlog} that the norm blows-up not only along a subsequence but pointwise for $t\rightarrow T^*$. This is due to the fact that the quantitative regularity under the a priori Type I control, see Result \ref{res.resquantest}, just involves the $L^3$ norm at the final time and not the whole $L^\infty_tL^3_x$ norm of the solution as the estimate \eqref{taomainest} of \cite{Tao19}.

\subsubsection*{\underline{Answer to Objective \ref{obj.quantreg}}}

\begin{resulta}[{\cite[Proposition 2.1, estimate (52)]{BP21cmp}}, localized quantitative regularity]\label{res.resquantest}
Assume that $U$ is a 
smooth solution with sufficient decay to the Navier-Stokes equations.\\ 
Assume in addition that $U$ satisfies the Type I bound
\begin{equation*}
A:=\|U\|_{L^{\infty}(-1,0;L^{3,\infty}(\mathbb{R}^3))}<\infty.
\end{equation*}
Then, letting
\begin{align*}\label{e.estdissipTypeI}
\begin{split}
t_*=\ &t_*(A,U(\cdot,0))\\
:=\ &-A^{-c}\exp\bigg(-\exp\exp(A^{1024})\int\limits_{B_{0}(\exp(A^{1023}))}|U(\cdot,0)|^3\bigg),
\end{split}
\end{align*}
the following quantitative boundedness holds
\begin{equation}\label{e.quantesttypeI}
\|U\|_{L^\infty(B_{0}(A^{c}\sqrt{-t_*})\times(t_*/2,0))}\lesssim \frac{A^{-c}}{\sqrt{-t_*}},
\end{equation}
with $c\in(0,\infty)$ a universal constant.
\end{resulta}

Notice that \eqref{e.quantesttypeI} involves the a priori Type I control in $L^\infty_tL^{3,\infty}_x$ and the control of the $L^3$ norm at the final time. Given the current state of the art, such assumptions are needed to prove the regularity; a mere Type I assumption is at this point not enough to beat the scaling except in the axisymmetric case and the self-similar case.\footnote{This is due to an additional scalar structure in those cases.} For more insights on the qualitative regularity in borderline endpoint critical spaces such as the Lorentz space $L^{3,\infty}$ or the Besov space $B^{-1+\frac3p}_{p,\infty}$, see \cite{AlbrittonBarkerBesov2018}.\footnote{These spaces are sometimes also referred to as `ultracritical spaces'. Unlike $L^3$, a function in these spaces can have a simultaneous presence in terms of its norm at an arbitrary amount of disjoint scales/frequencies, see \cite[Footnote 7]{Tobias2021}} 

Estimate \eqref{e.quantesttypeI} also implies a localized quantitative bound in the spirit of  \eqref{quantgenform} with $\mathcal B:=L^3(B_0(\exp(A^{1023})))$, $p:=\infty$, $\mathcal A:=L^{3,\infty}(\R^3)$, 
\begin{multline*}
\mathscr G(A,B)\simeq A^c\exp\big(\exp(\exp(A))B),\ \mbox{where}\  A:=\|U\|_{L^{\infty}(-1,0;L^{3,\infty}(\mathbb{R}^3))}\\ 
\mbox{and}\ B:=\int\limits_{B_0(\exp(A^{1023}))}|U(\cdot,0)|^3.
\end{multline*}

\subsubsection*{\underline{Answer to Objective \ref{goal.quantdissip}}}

\begin{resulta}[{\cite[Proposition 2.1, estimate (50)]{BP21cmp}}, quantitative estimates of dissipative/Kolmogorov scales]\label{res.reskol}
Assume that $U$ is a 
smooth solution with sufficient decay to the Navier-Stokes equations.\\ 
Assume in addition that $U$ satisfies the Type I bound
\begin{equation*}
A:=\|U\|_{L^{\infty}(-1,0;L^{3,\infty}(\mathbb{R}^3))}<\infty.
\end{equation*}
Then the threshold determining where the dissipative effects dominate the nonlinear effects is estimated as follows:
\begin{align}\label{e.estdissipTypeI}
\begin{split}
t_*=\ &t_*(A,U(\cdot,0))\\
:=\ &-A^{-c}\exp\bigg(-\exp\exp(A^{1024})\int\limits_{B_{0}(\exp(A^{1023}))}|U(\cdot,0)|^3\bigg).
\end{split}
\end{align}
\end{resulta}

It is not surprising that the definition of $t_*$ already appears in Result \ref{res.resquantest}, since Result \ref{res.resquantest} is (almost) an immediate consequence of Result \ref{res.reskol} and a regularity criteria (namely a local-in-space short-time regularity result). More precisely, $-t_*$ defined by \eqref{e.estdissipTypeI} is the time after which the following scale-invariant enstrophy\footnote{Here, as is usual, $\omega=\nabla\times U$.}
\begin{equation}\label{vortconcintro}
\mathscr S(A;t):=(-t)^{\frac{1}{2}}\int\limits_{B_{0}\big(A^c(-t)^{\frac 12}\big)}|\omega(x,t)|^2 dx
\end{equation}
is smaller than a given small number $\ep(A)=A^{-c}$.

Estimate \eqref{e.estdissipTypeI} also implies a quantitative bound like \eqref{goal.quantdissip} with $\mathcal B:=L^3(B_0(\exp(A^{1023})))$, $p:=\infty$, $\mathcal A:=L^{3,\infty}(\R^3)$, 
\begin{multline*}
\mathscr H(A,B)\simeq A^{-c}\exp\big(-\exp\exp(A^{1024})B\big),\ \mbox{where}\  A:=\|U\|_{L^{\infty}(-1,0;L^{3,\infty}(\mathbb{R}^3))}\\
\mbox{and}\ B:=\int\limits_{B_0(\exp(A^{1023}))}|U(\cdot,0)|^3.
\end{multline*}

\subsection{Quantitative estimates of dissipative scales in the Type I case}
\label{sec.typeIestdissscales}

The general scheme to get the estimate \eqref{e.estdissipTypeI} is a physical space parallel to Tao's Fourier based strategy described in Section \ref{sec.taosstrat}. Assume that $\mathscr S(A;t)$ defined by \eqref{vortconcintro} concentrates, i.e. is not small, for a certain time $t$. Then, the goal is to find an upper bound $t_*$ for $t$. More precisely, assume that 
\begin{equation}\label{e.initconcTypeI}
\mathscr S(A;t):=(-t)^{\frac{1}{2}}\int\limits_{B_{0}\big(A^c(-t)^{\frac 12}\big)}|\omega(x,t)|^2 dx>A^{-c},
\end{equation}
for a $t\in (-1,0)$ is not too close to $-1$. We call this the `initial concentration'.

\begin{figure}[h]
\begin{center}
\includegraphics[scale=.65]{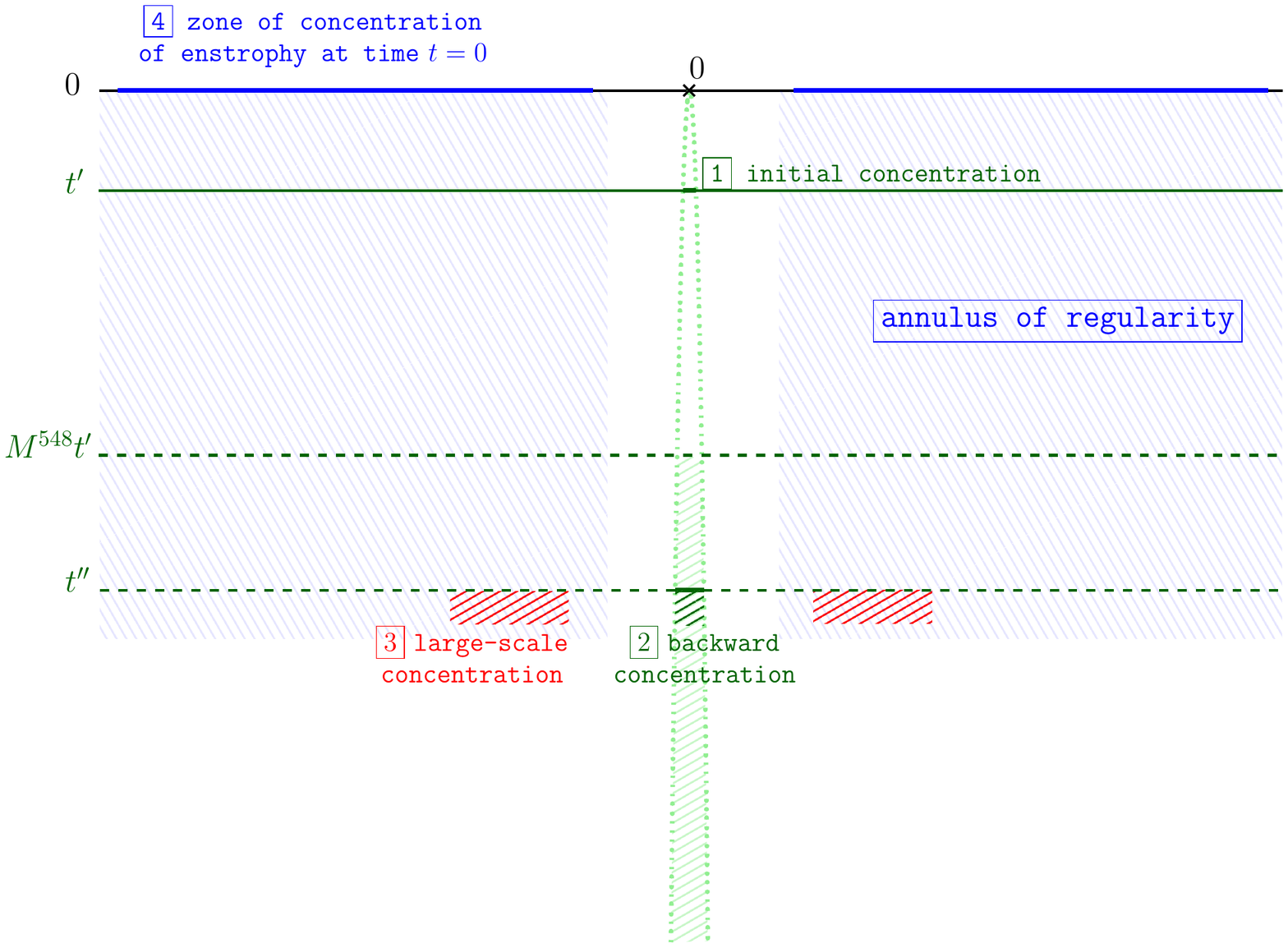}
\caption[Quantitative regularity via concentration: propagation of concentration]{Quantitative regularity in the Type I case via concentration of a scale-invariant enstrophy: \ref{TypeIstep1} backward, large-scale and forward propagation of enstrophy concentration}
\label{fig.c7BPprop}
\end{center}
\end{figure}

\begin{figure}[h]
\begin{center}
\includegraphics[scale=.65]{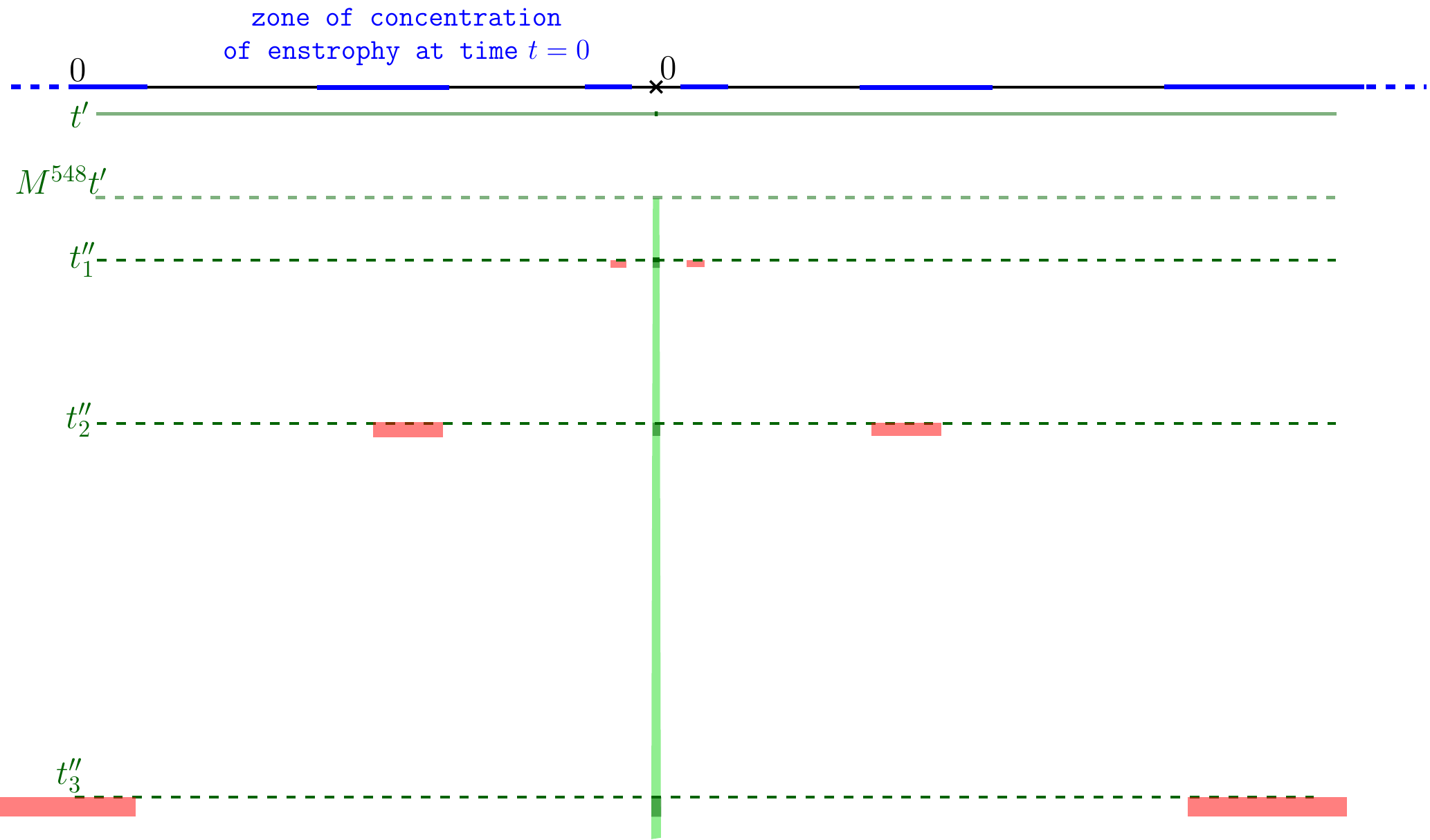}
\caption[Quantitative regularity via concentration: summation of scales]{Quantitative regularity in the Type I case via concentration of a scale-invariant enstrophy: \ref{TypeIstep2} summation of scales}
\label{fig.c7BP}
\end{center}
\end{figure}

\begin{enumerate}[label=(Step-\arabic*)]
\item \label{TypeIstep1}\textbf{Propagation of concentration}\\
For this step we refer to Figure \ref{fig.c7BP}. The idea is to transfer the initial concentration \eqref{e.initconcTypeI} in space and time in order to get a lower bound on the localized $L^3$ norm at time zero. That propagation of concentration relies on:
\begin{enumerate}[label=(\roman*)]
\item \emph{Backward propagation of concentration}\\
For all $t'\in (-1,t)$ such that $-t'$ is 
not too close to $-t$, we have
\begin{equation}\label{backwardvortconcintro}
(-t')^{\frac{1}{2}}\int\limits_{B_{0}\big(A^c(-t')^{\frac 12}\big)}|\omega(x,t')|^2 dx>A^{-c}.
\end{equation}
\item \emph{Large-scale and forward-in-time propagation of concentration}\\
It is shown that for certain admissible time scales $S$,\footnote{These admissible time scales are related to where one has backward concentration, see \cite[equation (106)-(107)]{BP21cmp}.} one has the concentration of the $L^3$ norm on the annulus $\{S^\frac12\leq|\cdot|\leq \exp(A^c)S^\frac12\}$, i.e.
\begin{equation}\label{L3lowerintro}
\int\limits_{S^\frac12\leq|\cdot|\leq \exp(A^c)S^\frac12} |U(x,0)|^3 dx\geq \exp(-\exp(A^{c})).
\end{equation}
The role of the Type I bound is to show that the solution $U$ obeys good quantitative estimates in certain space-time regions, epochs of quantitative regularity and annuli of quantitative regularity, which is needed to apply the Carleman inequalities to the vorticity equation, see \cite{ESS2003}, see  \cite[Proposition 4.2 and Proposition 4.3]{Tao19}. 
\end{enumerate}
\item \label{TypeIstep2}\textbf{Summation of scales and coercivity of the standing assumption}\\
We refer to Figure \ref{fig.c7BP} for this step. 
Summing \eqref{L3lowerintro} over all permissible disjoint annuli finally gives us the desired lower bound for $-t'$. We note that the localized $L^{3}$ norm of $U$ at time $0$ plays a distinct role to that of the Type I condition described in the previous step. Its sole purpose is to bound the number of permissible disjoint annuli that can be summed. This concludes the proof of the single exponential upper bound in ${\displaystyle \int\limits_{B_{0}(\exp(A^{1023}))}|U(\cdot,0)|^3}$ for $t$, see \eqref{e.estdissipTypeI}.
\end{enumerate}

Notice that quantitative local-in-space smoothing, see Section \ref{sec.locinspace} above, is a fundamental tool to achieve the backward propagation of enstrophy concentration in \ref{TypeIstep1} above as well as to go from Result \ref{res.reskol} to Result \ref{res.resquantest}.

\subsection{A comparison with Tao's strategy}

We outline here the two main differences between: on the one hand Tao's strategy \cite{Tao19} for the quantitative regularity in the case of a priori control of the solution in the borderline critical space $L^\infty_tL^3_x$ and on the other hand the strategy of \cite{BP21cmp} for the quantitative regularity in the case of a priori control of the solution in the borderline endpoint critical space $L^\infty_tL^{3,\infty}_x$ with additional control of the $L^3$ norm of $U(\cdot,0)$. There are two main aspects, see below. The first aspect is the decisive difference that also explains the second point. Let us state these two differences for the sake of clarity.

\subsubsection*{\underline{Fourier space vs. physical space}} 
We already underlined this aspect above. Tao works with the frequency bubbles of concentration \eqref{e.freqbubbles}, which are scale-invariant quantities defined in Fourier space. Therefore these quantities involve the solution $U(\cdot,t)$ on the whole-space $\R^3$. On the contrary, the analysis of \cite{BP21cmp} relies on the localized scale-invariant enstrophies \eqref{vortconcintro}. Those quantities are defined in physical space and involve the solution $U$ only on the ball $B_{0}(A^c(-t)^{\frac 12})$. Therefore, one can obtain localized results, such as the blow-up rate \eqref{L3localisedlog} of Result \ref{res.bluprateTypeI}.

\subsubsection*{\underline{Global scale-critical vs. criticality along a sequence of times/at a given time}} 
This aspect is a consequence of the point made above about Fourier vs. physical space approach. In Tao's work, there is a step, see `Transfer of concentration in Fourier space to physical space' in Section \ref{sec.taosstrat}, that consists in transferring the concentration of frequency bubbles to concentration of a scale-invariant enstrophy, see \cite[equation (5.6)]{Tao19}. That localized $L^2$ norm of the vorticity is needed to work with quantitative Carleman inequalities in order to propagate the concentration at large scales and forward in time. The process of transfer is based on the fact that the solution is bounded in $L^\infty(-1,0;L^3(\R^3))$. In the case of the strategy described in Section \ref{sec.typeIestdissscales} above in the Type I case, this step is not needed.

This remark has several implications:
\begin{enumerate}[label=(\arabic*)]
\item The scheme developed in the Type I case, see Section \ref{sec.typeIestdissscales}, enables to handle the case when a global scale-critical standing assumption is lacking, which is the case in \cite[Theorem B]{BP21cmp} that quantifies Seregin's 2012 criteria \cite{seregin2012}.
\item In a related vein, we get regularity at time $t=0$ for large $L^3$ data at $t=-1$ if the profile at time $t=0$ is quantitatively small; see \cite[Theorem 4.1 (i)]{AlbrittonBarkerBesov2018} for a qualitative statement and \cite[Proposition 4.4]{BP21cmp} for a quantitative statement.
\end{enumerate}

\section{Some further developments}
\label{sec.further}

There are many interesting research developments in the wake of the two papers \cite{Tao19,BP21cmp}. Let us review some of them.

Very recently, the paper \cite{barker22loc} was able to give a localized version of Tao's blow-up rate \eqref{e.quantbluprateTao}. For a smooth suitable weak solution $U$ to the Navier-Stokes equations on $B(0,4)\times (0,T^*)$, that possesses a singular point $(x_0,T^*)\in B(0,4)\times \{T^*\}$, then for all $\delta>0$ sufficiently small
$$\limsup_{t\uparrow T^*} \frac{\|U(\cdot,t)\|_{L^{3}(B(x_0,\delta))}}{\Big(\log\log\log\Big(\frac{1}{(T_*-t)^{\frac{1}{4}}}\Big)\Big)^{\frac{1}{1129}}}=\infty.$$
In that sense, this bound is a true quantification of the qualitative Escauriaza, Seregin and \v{S}ver\'{a}k \cite{ESS2003} criteria, which is also localized but qualitative. The main difficulty with the localization is related to showing that the local solution possesses quantitative annuli and epochs of regularity, which is required in \ref{step.step1Tao} \ref{step.step1Taoii}.

The quantitative regularity for solutions $U\in L^\infty_tL^d_x$ to the Navier-Stokes equations in higher dimensions $d\geq 4$ was handled by Palasek in \cite{Stan2}. This work gives an effective quantification of the qualitative result by Dong and Du \cite{DongDu}. For the blow-up rate, one pays the price of an additional logarithm compared to the result in dimension three \eqref{e.quantbluprateTao}.

Palasek \cite{Stan1} was also able to improve upon the triple logarithmic rate \eqref{e.quantbluprateTao} obtained by Tao in \cite{Tao19}: in the case of axisymmetric solutions for instance, the triple logarithm is replaced by a double logarithm. Without any symmetry assumption on the solution, a similar improvement can be obtained by replacing the $L^3$ norm by the scale-invariant norm $\|r^{1-\frac3q}U\|_{L^\infty_tL^q_x}$ for $q\in (3,\infty)$ and $r:=|x_h|$.

Finally, let us mention a related research line that aims at quantifying the regularity of axisymmetric solutions satisfying a critical or slightly supercritical Type I a priori bound. In the wake of the result of Pan \cite{pan2016regularity}, De Giorgi methods were intensively used to improve upon the regularity beyond the Type I case (see \cite{chen2008lower,chen2009lower,koch2009liouville,lei2011liouville}) by slightly breaking the a priori scale-invariant assumption. This research was carried by Seregin \cite{Seregin22,Seregin22bis} and Chen, Tsai and Zhang \cite{chen2022remarks}. In this last paper, a double-logarithmic quantitative blow-up rate for the $\dot B^{-1}_{\infty,\infty}(\R^3)$ norm of $U$ is obtained; see \cite[Theorem 1.4]{chen2022remarks}. Using Harnack inequalities instead of the Carleman inequalities used by Tao \cite{Tao19}, O{\.z}a{\'n}ski and Palasek \cite{OP22} recover the blow-up rate of the $L^{3,\infty}(\R^3)$ norm of $U$; see \cite[Corollary 1.2]{OP22}. In addition, they obtain a quantitative bound of the form \eqref{quantgenform} in terms of the $L^\infty_tL^{3,\infty}_x$ norm of the solution only; see \cite[Theorem 1.1]{OP22}. For other developments, see \cite{LeiRen}.

\section{Mild criticality breaking and a conjecture of Tao}
\label{sec.mildbreak}

Recently, in \cite[Remark 1.6]{Tao19}, Tao conjectured that if a solution first loses smoothness at time $T^*>0$, then the Orlicz norm $\|U(\cdot,t)\|_{L^{3}(\log\log\log L)^{-c}(\mathbb{R}^3)} $ must blow-up as $t$ tends to $T^*$. Result \ref{res.orliczblowup} provides a positive answer to Tao's conjecture, at the cost of one additional logarithm. 

\begin{resulta}[blow-up of slightly supercritical Orlicz norms; {\cite[Theorem 2]{BP21jmfm}}]
There exists a universal constant $\theta\in(0,1)$ such that the following holds.\\
Let $U$ be a Leray-Hopf solution to the Navier-Stokes equations on $\R^3\times(0,\infty)$ with initial data $U_0\in L^2(\R^3)\cap L^4(\R^3)$. Assume that $U$ first blows-up at $T^*\in(0,\infty)$. 
Then 
\begin{equation}\label{e.orliczblowup}
\lim\sup_{t\uparrow T^*}\int\limits_{\R^3}\frac{|U(x,t)|^3}{\Big(\log\log\log\big((\log(e^{e^{3e^{e}}}+|U(x,t)|))^\frac13\big)\Big)^\theta}dx=\infty.
\end{equation}
\label{res.orliczblowup}
\end{resulta}

As far as we know, Result \ref{res.orliczblowup} is the first result of this type for the Navier-Stokes equations concerned with slight criticality breaking in borderline spaces. Previously, it was shown for non-borderline spaces by Chan and Vasseur \cite{CV07} that if $U$ is a Leray-Hopf solution satisfying $$\int\limits_{0}^{\infty}\int\limits_{\mathbb{R}^3}\frac{|U|^5}{\log(1+|U|)}dxdt<\infty $$ then $U$ is smooth on $\mathbb{R}^3\times (0,\infty)$. Subsequent improvements were obtained in \cite{LZ13} and \cite{BV11}; see also \cite{MontgomerySmith05-mildbreak}. Let us mention that the techniques used in these papers cannot be used to treat the borderline case considered in Result \ref{res.orliczblowup}. %We also mention the paper of Chan and Yoneda \cite{ChanYoneda12} that mixes the boundedness of a strongly supercritical norm $L^\infty_tL^{\alpha,\infty}_x$ for $\alpha\in(2.343,3)$ with a geometrical information involving $\nabla\cdot(U/|U|)$.

\subsubsection*{\underline{A new method for transferring subcriticality of the data forward in time}}

The method for proving Result \ref{res.orliczblowup} relies on the following lemma and on a careful tuning of the parameters (estimating the $L^{3-\mu}$ norm for a well-chosen parameter $\mu$). Lemma \ref{lem.bulut} is directly inspired by the recent result of Bulut \cite{Bulut20} for a nonlinear supercritical defocusing Schr\"{o}dinger equation. 

\begin{resulta}[`mild criticality breaking'; {\cite[Theorem 1]{BP21jmfm}}]\label{lem.bulut}
For all $M,\, A\in[1,\infty)$ sufficiently large, 
there exists $\delta(M,A)\in(0,\frac12]$ such that the following holds. Let $U$ be a suitable weak Leray-Hopf solution to the Navier-Stokes equations on $\R^3\times(0,\infty)$ with initial data $U_0\in L^2(\R^3)\cap L^4(\R^3)$.\\
Assume that 
\begin{equation*}
\|U_0\|_{L^2},\ \|U_0\|_{L^4}\leq M,
\end{equation*}
and that
\begin{equation}\label{slightlysupercrithypo}
\|U\|_{L^\infty(0,\infty;L^{3-\delta(M,A)}(\R^3))}\leq A.
\end{equation}
Then, the above assumptions imply that $U$ is smooth on $\R^3\times (0,\infty)$. 
Moreover, there is an explicit formula for $\delta(M,A)$, see \cite[equation (26)]{BP21jmfm}, and $\delta(M,A)\rightarrow 0$ when $M\rightarrow\infty$ or $A\rightarrow\infty$.
\end{resulta}

We call Result \ref{lem.bulut} a `mild breaking of the criticality', or a `mild supercriti\-cal regularity criteria' as opposed to strong criticality breaking results obtained for instance in the axisymmetric case \cite{pan2016regularity,Seregin22,Seregin22bis,chen2022remarks}, see Section \ref{sec.further}. Indeed, the supercritical space $L^{\infty}_tL^{3-\delta(M,A)}$ in which we break the scaling depends on the size $A$ of the solution in this supercritical space via $\delta(M,A)$. In other words this can be considered as a non-effective regularity criteria, hence the term `mild'. Moreover, given a solution $U$, assume that you knew all the $L^\infty_tL^{3-\delta}_x$ norms for $\delta\rightarrow 0$. Then the question whether Lemma \ref{lem.bulut} applies to $U$ or not becomes a question about how fast 
\begin{equation*}
\|U\|_{L^\infty(0,\infty;L^{3-\delta}(\R^3))}
\end{equation*} 
grows when $\delta\rightarrow 0$. 
Of course one would have regularity if the solution was a priori bounded in the critical space $L^\infty_tL^3_x$. The result shows that with $L^4$ initial data one can relax the exponent $3$ to a slightly supercritical $3-\delta(M,A)$. Let us also remark that the initial condition $U_0\in L^4(\R^3)$ can be replaced by any subcritical initial condition $U_0\in L^{3+}(\R^3)$.

The main idea of the proof of Lemma \ref{lem.bulut} is to transfer subcritical information from the initial time forward in time. In a nutshell, subcritical energy estimates are combined with quantitative regularity estimates as obtained by Tao \cite{Tao19}. Hence the growth of the subcritical norm along the evolution can be estimated.\footnote{Result \ref{lem.bulut} can be abstractly quantified using persistence of singularities, see \cite[Introduction]{BP21jmfm}.} 
In that perspective the main objective is the following.

\begin{goala}
Prove that there exists $\delta(M,A)\in(0,\frac12]$ and $K(M,A)\in[1,\infty)$ such that for all $U_0\in L^2(\R^3)\cap L^4(\R^3)$ and any suitable weak Leray-Hopf solution associated to the initial data $U_0$, if 
\begin{equation*}
\|U_0\|_{L^2},\ \|U_0\|_{L^4}\leq M,
\end{equation*}
and
\begin{equation*}
\|U\|_{L^\infty(0,\infty;L^{3-\delta(M,A)}(\R^3))}\leq A,
\end{equation*}
then 
\begin{equation}\label{e.quantLinftyL4}
\|U\|_{L^\infty(0,T;L^4(\R^3))}\leq K(M,A).
\end{equation}
for any $T>0$.
\end{goala}
This then obviously implies the result stated in Result \ref{lem.bulut}. The crucial point is that $K(M,A)$ is uniform in time. 

Let us emphasize that the only a priori globally controlled quantity is a supercritical $L^\infty_tL^{3-}$ norm. We are not aware of any regularity mechanism enabling to break the critically barrier based on the sole knowledge of such a supercritical bound. 
Therefore, the idea, following Bulut \cite{Bulut20}, is to transfer the subcritical information coming from the initial data $U_0\in L^4(\R^3)$ to arbitrarily large times by using three ingredients:
\begin{enumerate}[label=(\arabic*)]
\item the control of the critical $L^\infty_tL^3_x$ norm via interpolation between the supercritical norm $L^\infty_tL^{3-\delta(M,A)}_x$ and the subcritical $L^\infty_tL^4_x$ norm
\begin{align*}
\begin{split}\label{e.estLinftyL3}
\|U\|_{L^\infty(0,T;L^3(\R^3))}\leq\ &\|U\|_{L^\infty(0,T;L^{3-\delta}(\R^3))}^\frac{3-\delta}{3+3\delta}\|U\|_{L^\infty(0,T;L^4(\R^3))}^\frac{4\delta}{3+3\delta}\\
\leq\ &A^\frac{3-\delta}{3+3\delta}K^\frac{4\delta}{3+3\delta}\, ;
\end{split}
\end{align*}
\item the quantitative control of the critical non borderline $L^5_{t,x}$ norm (see \cite[Proposition 3]{BP21jmfm}) in terms of the critical norm $\|U\|_{L^\infty(0,\infty;L^{3}(\R^3))}$, and the supercritical $L^2$ and subcritical $L^4$ norms of the initial data $U_0$
\begin{equation*}\label{e.quantestinterpol}
\|U\|_{L^5(0,T;L^5(\R^3))}\leq C(M)\exp\exp\exp\big(C_{univ}\big(A^\frac{3-\delta}{3+3\delta}K^\frac{4\delta}{3+3\delta}\big)^c\big);
\end{equation*} 
this hinges on the quantitative bounds on solutions belonging to the critical space $L^{\infty}_{t}L^{3}_{x}$, which were established by Tao in \cite{Tao19}, see \eqref{taomainest}; this step enables the slicing of the interval $(0,T)$ into a $T$-independent number $m$ of disjoint epochs $I_j=(t_j,t_{j+1})$,
\begin{align*}
\ep^5m=\sum_{j=1}^m\|U\|_{L^5(I_{j};L^5(\R^3))}^5\leq \ &\|U\|_{L^5(0,T;L^5(\R^3))}^5\\
\leq\ &C(M)\exp\exp\exp\Big(\big(A^\frac{3-\delta}{3+3\delta}K^\frac{4\delta}{3+3\delta}\big)^c\Big)\, ;
\end{align*}
\item an $L^4$ energy estimate \cite[Proposition 4]{BP21jmfm} under the $L^5_{t,x}$ control of $U$, which enables the transfer the subcritical information from time $t_j$ to $t_{j+1}$
\begin{align*}
\mathscr E_{4,t_{j+1}}\leq\ & \|U(\cdot,t_j)\|_{L^4(\R^3)}^4+C\|U\|_{L^5(\R^3\times I_j)}\mathcal E_{4,t_{j+1}}\\
\leq\ &\|U(\cdot,t_{j+1})\|_{L^4(\R^3)}^4+C\ep\mathscr E_{4,t_{j+1}},
\end{align*} 
where $\mathscr E_{4,t_{j+1}}$ is the $L^4$ energy, see \cite[equation (13)]{BP21jmfm}, and eventually to $T$
\begin{align*}
\|U\|_{L^\infty(0,T;L^4(\R^3))}^4=\max_{1\leq j\leq m+1}\{\|U\|_{L^\infty(I_j;L^4(\R^3))}^4\}\leq 64M^4 2^m.
\end{align*}
One then designs the number $K(M,A)$ and $\delta(M,A)$ to bound the right hand side above.
\end{enumerate}

\section{Summary of selected results}
\label{sec.summary}

Figure \ref{fig.graphconc} on page \pageref{fig.graphconc} summarizes some results about weak and strong concentration for the Navier-Stokes equations, as well as local-in-space smoothing. Figure \ref{fig.graphqq} on page \pageref{fig.graphqq} summarizes qualitative and quantitative regularity results for the Navier-Stokes equations.

\begin{figure}[h]
\vspace{1cm}
\begin{adjustwidth}{-1.5cm
}{}
\begin{tabular}{L{.5cm}|L{6.3cm}|L{6cm}|L{5.3cm}}
&\textbf{weak concentration}&\textbf{local-in-space smoothing}&\textbf{strong concentration}\\
\hline
\rotatebox[origin=c]{90}{\textbf{local energy solutions/Type I}}&&\shortstack[l]{Jia, {\v{S}}ver{\'a}k (2014) \cite{JS14} \emph{subcritical}\\Barker, Prange (2020) \cite{BP18}\\ \emph{critical, whole-space}\\ Kang, Miura, Tsai (2020) \cite{KMT18}\\ \emph{critical, whole-space}\\Kang, Miura, Tsai (2020) \cite{KangMiuraTsai20-concL2}\\ \emph{scaled energy, whole-space}\\ Albritton, Barker, Prange (2021) \cite{ABP21}\\ \emph{critical, half-space}}&\shortstack[l]{Barker, Prange (2020) \cite{BP18}\\ \emph{whole-space}\\ Albritton, Barker, Prange (2021) \cite{ABP21}\\ \emph{half-space}\\ Barker, Prange (2021) \cite{BP21cmp}\\ \emph{blow-up of}\\
$L^3\big(B_{0}(O((T^*-t)^{\frac12-}))\big)$\\ Kang, Miura, Tsai (2021) \cite{KangMiuraTsai20-concL2}\\
\emph{scale-invariant energy}}\\
\hline  
\rotatebox[origin=c]{90}{\textbf{beyond Type I}} & \shortstack[l]{Li, Ozawa, Wang (2018) \cite{LOW18}\\ Maekawa, Miura, Prange (2020) \cite{MMP17a}\\ Kang, Miura, Tsai (2021) \cite{KangMiuraTsai20-concL2}\\ \emph{scale-invariant energy}}&&
\end{tabular}
\end{adjustwidth}
\caption[Concentration]{Weak and strong concentration: a selection of results}
\label{fig.graphconc}
\end{figure}

\begin{figure}[h]
\vspace{1cm}
\begin{adjustwidth}{-2cm}{}
\begin{tabular}{L{1cm}|L{3cm}|L{3cm}|L{3cm}|L{4.5cm}|L{4cm}}
&\textbf{critical\tikzmark{e}}&\textbf{borderline critical}&\textbf{borderline endpoint critical}&\textbf{slightly supercritical}&\textbf{supercritical}\tikzmark{f}\\
&&$L^\infty_tL^3_x$&Type I&`log' breaking&$\sup_{k}\|U(\cdot,t_k)\|_{L^3}<\infty$\\
\hline  
\rotatebox[origin=c]{90}{\textbf{qualitative}}&Lady\v{z}enskaja-Prodi-Serrin (60's)&Escauriaza, Seregin, \v{S}ver\'{a}k (2003) \cite{ESS2003}, Kenig, Koch (2011) \cite{KKH1/2}, Gallagher, Koch, Planchon (2013) \cite{GKP}\ldots & &$\int\limits
_0^T\int\limits_{\R^3}\frac{|U|^5}{\log(1+|U|)}\, dxds<\infty$ Chan, Vasseur (2007) \cite{CV07}, Bjorland, Vasseur (2011) \cite{BV11}, Lei, Zhou (2013) \cite{LZ13}\ldots&Seregin (2012) \cite{seregin2012}, Albritton and Barker \cite[Theorem 4.1 (i)]{AlbrittonBarkerBesov2018}\\
&&&&Barker, Prange (2021) \cite{BP21jmfm} \emph{blow-up of an Orlicz norm}%\footnote{The proof of this result is quantitative in nature so we expect a corresponding blow-up rate.}
&\\
%\hline
%\rotatebox[origin=c]{90}{\textbf{abstract quantitative}}&&Seregin (2010's)&&\\
\hline
\rotatebox[origin=c]{90}{\textbf{explicit quantitative}}&&Tao (2019) \cite{Tao19}&Barker, Prange (2021) \cite{BP21cmp} \emph{localized blow-up of $L^3$}&&Barker, Prange (2021) \cite{BP21cmp} \emph{quantification of Seregin's 2012 result and of Albritton and Barker's 2019 result}
\end{tabular}
\begin{tikzpicture}[remember picture,overlay]
\draw[-latex] ([ shift={(0,6ex)}]pic cs:e) -- node[above] {increasing level of criticality} ([ shift={(0,6ex)}]pic cs:f);
\end{tikzpicture}
\end{adjustwidth}
\caption[Qualitative vs. quantitative regularity]{Qualitative vs. abstract and explicit quantitative regularity: a selection of results}
\label{fig.graphqq}
\end{figure}            

\subsection*{Acknowledgement}
Both authors thank the Institute of Advanced Studies of Cergy Paris University for their hospitality. CP is partially supported by the Agence Nationale de la Recherche,
project BORDS, grant ANR-16-CE40-0027-01, project SINGFLOWS, grant ANR-
18-CE40-0027-01, project CRISIS, grant ANR-20-CE40-0020-01, by the CY
Initiative of Excellence, project CYNA (CY Nonlinear Analysis) and project CYFI (CYngular Fluids and Interfaces).

\small 
\bibliographystyle{abbrv}
\bibliography{concentration.bib}

\begin{thebibliography}{10}

\bibitem{AlbrittonBarkerBesov2018}
D.~Albritton and T.~Barker.
\newblock Global weak besov solutions of the {N}avier-{S}tokes equations and
  applications.
\newblock {\em Archive for Rational Mechanics and Analysis}, 232(1):197--263,
  2019.

\bibitem{ABP21}
D.~{Albritton}, T.~{Barker}, and C.~{Prange}.
\newblock {Localized smoothing and concentration for the {N}avier-{S}tokes
  equations in the half space}.
\newblock {\em arXiv e-prints}, page arXiv:2112.10705, Dec. 2021.

\bibitem{AC10-concentration}
M.~Arnold and W.~Craig.
\newblock On the size of the {Navier}-{Stokes} singular set.
\newblock {\em Discrete Contin. Dyn. Syst.}, 28(3):1165--1178, 2010.

\bibitem{Tobias2021}
T.~{Barker}.
\newblock {Higher integrability and the number of singular points for the
  Navier-Stokes equations with a scale-invariant bound}.
\newblock {\em arXiv e-prints}, page arXiv:2111.14776, Nov. 2021.

\bibitem{barker22loc}
T.~{Barker}.
\newblock {Localized quantitative estimates and potential blow-up rates for the
  Navier-Stokes equations}.
\newblock {\em arXiv e-prints}, page arXiv:2209.15627, Sept. 2022.

\bibitem{BP18}
T.~Barker and C.~Prange.
\newblock Localized {S}moothing for the {N}avier-{S}tokes {E}quations and
  {C}oncentration of {C}ritical {N}orms {N}ear {S}ingularities.
\newblock {\em Arch. Ration. Mech. Anal.}, 236(3):1487--1541, 2020.

\bibitem{BP20geom}
T.~Barker and C.~Prange.
\newblock Scale-invariant estimates and vorticity alignment for
  {N}avier-{S}tokes in the half-space with no-slip boundary conditions.
\newblock {\em Arch. Ration. Mech. Anal.}, 235(2):881--926, 2020.

\bibitem{BP21jmfm}
T.~Barker and C.~Prange.
\newblock Mild criticality breaking for the {N}avier-{S}tokes equations.
\newblock {\em J. Math. Fluid Mech.}, 23(3):Paper No. 66, 12, 2021.

\bibitem{BP21cmp}
T.~Barker and C.~Prange.
\newblock Quantitative regularity for the {N}avier-{S}tokes equations via
  spatial concentration.
\newblock {\em Comm. Math. Phys.}, 385(2):717--792, 2021.

\bibitem{BV11}
C.~Bjorland and A.~Vasseur.
\newblock Weak in space, log in time improvement of the
  {L}ady\v{z}enskaja-{P}rodi-{S}errin criteria.
\newblock {\em J. Math. Fluid Mech.}, 13(2):259--269, 2011.

\bibitem{Bourgain98}
J.~Bourgain.
\newblock Refinements of {S}trichartz' inequality and applications to
  {$2$}{D}-{NLS} with critical nonlinearity.
\newblock {\em Internat. Math. Res. Notices}, (5):253--283, 1998.

\bibitem{bradshaw2022local}
Z.~Bradshaw and T.-P. Tsai.
\newblock On the local pressure expansion for the {N}avier-{S}tokes equations.
\newblock {\em Journal of Mathematical Fluid Mechanics}, 24(1):1--32, 2022.

\bibitem{Bulut20}
A.~{Bulut}.
\newblock {Blow-up criteria below scaling for defocusing energy-supercritical
  NLS and quantitative global scattering bounds}.
\newblock {\em arXiv e-prints}, page arXiv:2001.05477, Jan. 2020.

\bibitem{CKN82}
L.~Caffarelli, R.~Kohn, and L.~Nirenberg.
\newblock Partial regularity of suitable weak solutions of the
  {N}avier-{S}tokes equations.
\newblock {\em Comm. Pure Appl. Math.}, 35(6):771--831, 1982.

\bibitem{camliyurt2022scattering}
G.~Camliyurt and C.~E. Kenig.
\newblock Scattering for focusing supercritical wave equations in odd
  dimensions.
\newblock {\em arXiv preprint arXiv:2201.04710}, 2022.

\bibitem{cannone1997generalization}
M.~Cannone.
\newblock A generalization of a theorem by {K}ato on {N}avier-{S}tokes
  equations.
\newblock {\em Revista matem{\'a}tica iberoamericana}, 13(3):515--541, 1997.

\bibitem{cazenave1989introduction}
T.~Cazenave.
\newblock {\em An introduction to nonlinear {S}chr{\"o}dinger equations},
  volume~22.
\newblock Universidade federal do Rio de Janeiro, Centro de ci{\^e}ncias
  matem{\'a}ticas e da~?, 1989.

\bibitem{CW17}
D.~Chae and J.~Wolf.
\newblock Removing discretely self-similar singularities for the 3{D}
  {N}avier-{S}tokes equations.
\newblock {\em Comm. Partial Differential Equations}, 42(9):1359--1374, 2017.

\bibitem{chae2020energy}
D.~Chae and J.~Wolf.
\newblock Energy concentrations and {T}ype {I} blow-up for the {3D} {E}uler
  equations.
\newblock {\em Communications in Mathematical Physics}, 376(2):1627--1669,
  2020.

\bibitem{CV07}
C.~H. Chan and A.~Vasseur.
\newblock Log improvement of the {P}rodi-{S}errin criteria for
  {N}avier-{S}tokes equations.
\newblock {\em Methods Appl. Anal.}, 14(2):197--212, 2007.

\bibitem{chen2009lower}
C.-C. Chen, R.~M. Strain, T.-P. Tsai, and H.-T. Yau.
\newblock Lower bounds on the blow-up rate of the axisymmetric
  {N}avier-{S}tokes equations {II}.
\newblock {\em Communications in Partial Differential Equations},
  34(3):203--232, 2009.

\bibitem{chen2008lower}
C.-C. Chen, R.~M. Strain, H.-T. Yau, and T.-P. Tsai.
\newblock Lower bound on the blow-up rate of the axisymmetric {N}avier-{S}tokes
  equations.
\newblock {\em International Mathematics Research Notices}, 2008, 2008.

\bibitem{chen2022remarks}
H.~Chen, T.-P. Tsai, and T.~Zhang.
\newblock Remarks on local regularity of axisymmetric solutions to the {3D}
  {N}avier-{S}tokes equations.
\newblock {\em Communications in Partial Differential Equations}, pages 1--20,
  2022.

\bibitem{ChesDai19}
A.~Cheskidov and M.~Dai.
\newblock Kolmogorov's dissipation number and the number of degrees of freedom
  for the 3d {N}avier-{S}tokes equations.
\newblock {\em Proc. R. Soc. Edinb., Sect. A, Math.}, 149(2):429--446, 2019.

\bibitem{cheskidov2010regularity}
A.~Cheskidov and R.~Shvydkoy.
\newblock The regularity of weak solutions of the 3d {N}avier-{S}tokes
  equations in {$B^{-1}_{\infty,\infty}$}.
\newblock {\em Archive for rational mechanics and analysis}, 195(1):159--169,
  2010.

\bibitem{CS14}
A.~Cheskidov and R.~Shvydkoy.
\newblock A unified approach to regularity problems for the 3{D}
  {N}avier-{S}tokes and {E}uler equations: the use of {K}olmogorov's
  dissipation range.
\newblock {\em J. Math. Fluid Mech.}, 16(2):263--273, 2014.

\bibitem{DongDu}
H.~Dong and D.~Du.
\newblock The {N}avier-{S}tokes equations in the critical {L}ebesgue space.
\newblock {\em Comm. Math. Phys.}, 292(3):811--827, 2009.

\bibitem{DY18}
T.~Duyckaerts and J.~Yang.
\newblock Blow-up of a critical {S}obolev norm for energy-subcritical and
  energy-supercritical wave equations.
\newblock {\em Anal. PDE}, 11(4):983--1028, 2018.

\bibitem{ESS2003}
L.~Escauriaza, G.~A. Seregin, and V.~\v{S}ver\'{a}k.
\newblock {$L_{3,\infty}$}-solutions of {N}avier-{S}tokes equations and
  backward uniqueness.
\newblock {\em Uspekhi Mat. Nauk}, 58(2(350)):3--44, 2003.

\bibitem{GKP}
I.~Gallagher, G.~S. Koch, and F.~Planchon.
\newblock Blow-up of critical {B}esov norms at a potential {N}avier-{S}tokes
  singularity.
\newblock {\em Comm. Math. Phys.}, 343(1):39--82, 2016.

\bibitem{GrujicXu2019-dynrestr}
Z.~Gruji\'{c} and L.~Xu.
\newblock A regularity criterion for {3D} {NSE} in dynamicaly restricted local
  {M}orrey spaces.
\newblock {\em Applicable Analysis}, 0(0):1--15, 2021.

\bibitem{HK06}
T.~Hmidi and S.~Keraani.
\newblock Remarks on the blowup for the {$L^2$}-critical nonlinear
  {S}chr\"{o}dinger equations.
\newblock {\em SIAM J. Math. Anal.}, 38(4):1035--1047, 2006.

\bibitem{HR07}
J.~Holmer and S.~Roudenko.
\newblock On blow-up solutions to the 3{D} cubic nonlinear {S}chr\"{o}dinger
  equation.
\newblock {\em Appl. Math. Res. Express. AMRX}, (1):Art. ID abm004, 31, 2007.

\bibitem{hou2022potentially}
T.~Y. Hou.
\newblock Potentially singular behavior of the {3D} {N}avier--{S}tokes
  equations.
\newblock {\em Foundations of Computational Mathematics}, pages 1--49, 2022.

\bibitem{JS14}
H.~Jia and V.~\v{S}ver\'{a}k.
\newblock Local-in-space estimates near initial time for weak solutions of the
  {N}avier-{S}tokes equations and forward self-similar solutions.
\newblock {\em Invent. Math.}, 196(1):233--265, 2014.

\bibitem{Kang05}
K.~Kang.
\newblock Unbounded normal derivative for the {S}tokes system near boundary.
\newblock {\em Math. Ann.}, 331(1):87--109, 2005.

\bibitem{kang2021finite}
K.~Kang, B.~Lai, C.-C. Lai, and T.-P. Tsai.
\newblock Finite energy {N}avier-{S}tokes flows with unbounded gradients
  induced by localized flux in the half-space.
\newblock {\em arXiv preprint arXiv:2107.00810}, 2021.

\bibitem{KMT18}
K.~Kang, H.~Miura, and T.-P. Tsai.
\newblock {Short Time Regularity of {N}avier-{S}tokes Flows with Locally
  ${L}^3$ Initial Data and Applications}.
\newblock {\em International Mathematics Research Notices}, 01 2020.
\newblock rnz327.

\bibitem{KangMiuraTsai20-concL2}
K.~Kang, H.~Miura, and T.-P. Tsai.
\newblock Regular sets and an {$\epsilon$}-regularity theorem in terms of
  initial data for the {N}avier-{S}tokes equations.
\newblock {\em Pure Appl. Anal.}, 3(3):567--594, 2021.

\bibitem{KKH1/2}
C.~E. Kenig and G.~S. Koch.
\newblock An alternative approach to regularity for the {N}avier-{S}tokes
  equations in critical spaces.
\newblock {\em Ann. Inst. H. Poincar\'{e} C Anal. Non Lin\'{e}aire},
  28(2):159--187, 2011.

\bibitem{KM10}
C.~E. Kenig and F.~Merle.
\newblock Scattering for {$\dot H^{1/2}$} bounded solutions to the cubic,
  defocusing {NLS} in 3 dimensions.
\newblock {\em Trans. Amer. Math. Soc.}, 362(4):1937--1962, 2010.

\bibitem{KM11}
C.~E. Kenig and F.~Merle.
\newblock Nondispersive radial solutions to energy supercritical non-linear
  wave equations, with applications.
\newblock {\em Amer. J. Math.}, 133(4):1029--1065, 2011.

\bibitem{KPV00}
C.~E. Kenig, G.~Ponce, and L.~Vega.
\newblock On the concentration of blow up solutions for the generalized {K}d{V}
  equation critical in {$L^2$}.
\newblock In {\em Nonlinear wave equations ({P}rovidence, {RI}, 1998)}, volume
  263 of {\em Contemp. Math.}, pages 131--156. Amer. Math. Soc., Providence,
  RI, 2000.

\bibitem{KV11}
R.~Killip and M.~Visan.
\newblock The defocusing energy-supercritical nonlinear wave equation in three
  space dimensions.
\newblock {\em Trans. Amer. Math. Soc.}, 363(7):3893--3934, 2011.

\bibitem{koch2009liouville}
G.~Koch, N.~Nadirashvili, G.~Seregin, and V.~{\v{S}}ver{\'a}k.
\newblock Liouville theorems for the {N}avier-{S}tokes equations and
  applications.
\newblock {\em Acta Mathematica}, 203(1):83--105, 2009.

\bibitem{kwon2021role}
H.~Kwon.
\newblock The role of the pressure in the regularity theory for the
  {N}avier-{S}tokes equations.
\newblock {\em arXiv preprint arXiv:2104.03160}, 2021.

\bibitem{LeiRen}
Z.~{Lei} and X.~{Ren}.
\newblock {Quantitative partial regularity of the Navier-Stokes equations and
  applications}.
\newblock {\em arXiv e-prints}, page arXiv:2210.01783, Oct. 2022.

\bibitem{lei2011liouville}
Z.~Lei and Q.~S. Zhang.
\newblock A {L}iouville theorem for the axially-symmetric {N}avier-{S}tokes
  equations.
\newblock {\em Journal of Functional Analysis}, 261(8):2323--2345, 2011.

\bibitem{LZ13}
Z.~Lei and Y.~Zhou.
\newblock Logarithmically improved criteria for {E}uler and {N}avier-{S}tokes
  equations.
\newblock {\em Commun. Pure Appl. Anal.}, 12(6):2715--2719, 2013.

\bibitem{Leray}
J.~Leray.
\newblock Sur le mouvement d'un liquide visqueux emplissant l'espace.
\newblock {\em Acta Math.}, 63(1):193--248, 1934.

\bibitem{LOW18}
K.~Li, T.~Ozawa, and B.~Wang.
\newblock Dynamical behavior for the solutions of the {N}avier-{S}tokes
  equation.
\newblock {\em Commun. Pure Appl. Anal.}, 17(4):1511--1560, 2018.

\bibitem{Lin}
F.~Lin.
\newblock A new proof of the {C}affarelli-{K}ohn-{N}irenberg theorem.
\newblock {\em Comm. Pure Appl. Math.}, 51(3):241--257, 1998.

\bibitem{MMP17a}
Y.~{Maekawa}, H.~{Miura}, and C.~{Prange}.
\newblock {Estimates for the {N}avier-{S}tokes equations in the half-space for
  non localized data}.
\newblock {\em ArXiv e-prints}, Nov. 2017.

\bibitem{MMP17b}
Y.~Maekawa, H.~Miura, and C.~Prange.
\newblock Local energy weak solutions for the {N}avier-{S}tokes equations in
  the half-space.
\newblock {\em Comm. Math. Phys.}, 367(2):517--580, 2019.

\bibitem{maekawa2006}
Y.~Maekawa and Y.~Terasawa.
\newblock The {N}avier-{S}tokes equations with initial data in uniformly local
  $l^p$ spaces.
\newblock {\em Differential Integral Equations}, 19(4):369--400, 2006.

\bibitem{Merle92}
F.~Merle.
\newblock On uniqueness and continuation properties after blow-up time of
  self-similar solutions of nonlinear {S}chr\"{o}dinger equation with critical
  exponent and critical mass.
\newblock {\em Comm. Pure Appl. Math.}, 45(2):203--254, 1992.

\bibitem{Merle93}
F.~Merle.
\newblock Determination of blow-up solutions with minimal mass for nonlinear
  {S}chr\"{o}dinger equations with critical power.
\newblock {\em Duke Math. J.}, 69(2):427--454, 1993.

\bibitem{MR06}
F.~Merle and P.~Rapha{\"e}l.
\newblock Blow up of the critical norm for some radial ${L}^2$ super critical
  nonlinear {S}chr{\"o}dinger equations.
\newblock {\em American journal of mathematics}, 130(4):945--978, 2008.

\bibitem{MT90}
F.~Merle and Y.~Tsutsumi.
\newblock {$L^2$} concentration of blow-up solutions for the nonlinear
  {S}chr\"{o}dinger equation with critical power nonlinearity.
\newblock {\em J. Differential Equations}, 84(2):205--214, 1990.

\bibitem{MZ05}
F.~Merle and H.~Zaag.
\newblock Determination of the blow-up rate for a critical semilinear wave
  equation.
\newblock {\em Math. Ann.}, 331(2):395--416, 2005.

\bibitem{MT22}
H.~{Miura} and J.~{Takahashi}.
\newblock {Blow-up of the critical norm for a supercritical semilinear heat
  equation}.
\newblock {\em arXiv e-prints}, page arXiv:2206.10790, June 2022.

\bibitem{MizoguchiSouplet19}
N.~Mizoguchi and P.~Souplet.
\newblock Optimal condition for blow-up of the critical {$L^q$} norm for the
  semilinear heat equation.
\newblock {\em Adv. Math.}, 355:106763, 24, 2019.

\bibitem{MontgomerySmith05-mildbreak}
S.~Montgomery-Smith.
\newblock Conditions implying regularity of the three dimensional
  {N}avier-{S}tokes equation.
\newblock {\em Applications of Mathematics}, 50(5):451--464, 2005.

\bibitem{Nawa90}
H.~Nawa.
\newblock ``{M}ass concentration'' phenomenon for the nonlinear
  {S}chr\"{o}dinger equation with the critical power nonlinearity. {II}.
\newblock {\em Kodai Math. J.}, 13(3):333--348, 1990.

\bibitem{Nawa92}
H.~Nawa.
\newblock ``{M}ass concentration'' phenomenon for the nonlinear
  {S}chr\"{o}dinger equation with the critical power nonlinearity.
\newblock {\em Funkcial. Ekvac.}, 35(1):1--18, 1992.

\bibitem{NT98}
H.~Nawa and M.~Tsutsumi.
\newblock On blowup for the pseudo-conformally invariant nonlinear
  {S}chr\"{o}dinger equation. {II}.
\newblock {\em Comm. Pure Appl. Math.}, 51(4):373--383, 1998.

\bibitem{OP22}
W.~S. {O{\.z}a{\'n}ski} and S.~{Palasek}.
\newblock {Quantitative control of solutions to axisymmetric Navier-Stokes
  equations in terms of the weak $L^3$ norm}.
\newblock {\em arXiv e-prints}, page arXiv:2210.10030, Oct. 2022.

\bibitem{Stan1}
S.~Palasek.
\newblock Improved quantitative regularity for the {N}avier-{S}tokes equations
  in a scale of critical spaces.
\newblock {\em Archive for Rational Mechanics and Analysis},
  242(3):1479–1531, Sep 2021.

\bibitem{Stan2}
S.~Palasek.
\newblock A minimum critical blowup rate for the high-dimensional
  {N}avier-{S}tokes equations, 2021.

\bibitem{pan2016regularity}
X.~Pan.
\newblock Regularity of solutions to axisymmetric {N}avier-{S}tokes equations
  with a slightly supercritical condition.
\newblock {\em Journal of Differential Equations}, 260(12):8485--8529, 2016.

\bibitem{CP-hdr}
C.~Prange.
\newblock {\em {Concentration and quantitative regularity in homogenization and
  hydrodynamics}}.
\newblock Habilitation {\`a} diriger des recherches, {CY Cergy Paris
  Universit{\'e}}, Oct. 2022.
\newblock
  https://hal.archives-ouvertes.fr/tel-03814044/file/hdr-Prange-web.pdf.

\bibitem{rusin2011minimal}
W.~Rusin and V.~{\v{S}}ver{\'a}k.
\newblock Minimal initial data for potential {N}avier-{S}tokes singularities.
\newblock {\em Journal of Functional Analysis}, 260(3):879--891, 2011.

\bibitem{seregin2012}
G.~Seregin.
\newblock A certain necessary condition of potential blow up for
  {N}avier-{S}tokes equations.
\newblock {\em Communications in Mathematical Physics}, 312(3):833--845, 2012.

\bibitem{Seregin22}
G.~Seregin.
\newblock A note on local regularity of axisymmetric solutions to the
  {N}avier-{S}tokes equations.
\newblock {\em J. Math. Fluid Mech.}, 24(1):Paper No. 27, 13, 2022.

\bibitem{Seregin22bis}
G.~Seregin.
\newblock A slightly supercritical condition of regularity of axisymmetric
  solutions to the {N}avier-{S}tokes equations.
\newblock {\em J. Math. Fluid Mech.}, 24(1):Paper No.18, 17, 2022.

\bibitem{SSv10}
G.~Seregin and V.~{\v{S}}ver{\'a}k.
\newblock On a bounded shear flow in half-space.
\newblock {\em Zap. Nauchn. Sem. S.-Peterburg. Otdel. Mat. Inst. Steklov.
  (POMI)}, 385(Kraevye Zadachi Matematichesko\u{\i} Fiziki i Smezhnye Voprosy
  Teorii Funktsi\u{\i}. 41):200--205, 236, 2010.

\bibitem{Seregin2018}
G.~Seregin and V.~{\v{S}}ver{\'a}k.
\newblock Regularity criteria for {N}avier-{S}tokes solutions.
\newblock In Y.~Giga and A.~Novotn{\'y}, editors, {\em Handbook of Mathematical
  Analysis in Mechanics of Viscous Fluids}, pages 829--867. Springer
  International Publishing, Cham, 2018.

\bibitem{SerZajac}
G.~A. Seregin and W.~Zajaczkowski.
\newblock A sufficient condition of local regularity for the {N}avier-{S}tokes
  equations.
\newblock {\em Zap. Nauchn. Sem. S.-Peterburg. Otdel. Mat. Inst. Steklov.
  (POMI)}, 336(Kraev. Zadachi Mat. Fiz. i Smezh. Vopr. Teor. Funkts.
  37):46--54, 274, 2006.

\bibitem{Serrin1962}
J.~Serrin.
\newblock On the interior regularity of weak solutions of the {N}avier-{S}tokes
  equations.
\newblock {\em Archive for Rational Mechanics and Analysis}, 9(1):187--195, Jan
  1962.

\bibitem{sulem2007nonlinear}
C.~Sulem and P.-L. Sulem.
\newblock {\em The nonlinear {S}chr{\"o}dinger equation: self-focusing and wave
  collapse}, volume 139.
\newblock Springer Science \& Business Media, 2007.

\bibitem{T90}
S.~Takahashi.
\newblock On interior regularity criteria for weak solutions of the
  {N}avier-{S}tokes equations.
\newblock {\em Manuscripta Math.}, 69(3):237--254, 1990.

\bibitem{Tao06}
T.~Tao.
\newblock {\em Nonlinear dispersive equations}, volume 106 of {\em CBMS
  Regional Conference Series in Mathematics}.
\newblock Published for the Conference Board of the Mathematical Sciences,
  Washington, DC; by the American Mathematical Society, Providence, RI, 2006.
\newblock Local and global analysis.

\bibitem{tao2013localisation}
T.~Tao.
\newblock Localisation and compactness properties of the {N}avier-{S}tokes
  global regularity problem.
\newblock {\em Analysis \& PDE}, 6(1):25--107, 2013.

\bibitem{Tao19}
T.~Tao.
\newblock Quantitative bounds for critically bounded solutions to the
  {N}avier-{S}tokes equations.
\newblock In A.~Kechris, N.~Makarov, D.~Ramakrishnan, and X.~Zhu, editors, {\em
  Nine Mathematical Challenges: An Elucidation}, volume 104. American
  Mathematical Society, 2021.

\bibitem{Tsutsumi90}
Y.~Tsutsumi.
\newblock Rate of {$L^2$} concentration of blow-up solutions for the nonlinear
  {S}chr\"{o}dinger equation with critical power.
\newblock {\em Nonlinear Anal.}, 15(8):719--724, 1990.

\bibitem{weinstein1989nonlinear}
M.~I. Weinstein.
\newblock The nonlinear {S}chrodinger equation-singularity formation, stability
  and dispersion, the connection between infinite and finite dimensional
  dynamical systems.
\newblock {\em Contemp. Math.}, 99:213--232, 1989.

\end{thebibliography}

\end{document}